\documentclass{amsart}

\usepackage{color}
\usepackage{amssymb, amsmath, amsthm}
\usepackage{graphicx}
\usepackage{url, hyperref}
\usepackage{fullpage}
\usepackage{subfig}

\theoremstyle{plain}
\newtheorem{definition}{Definition}

\theoremstyle{plain}
\newtheorem{theorem}{Theorem}

\theoremstyle{plain}
\newtheorem{prop}{Proposition}

\theoremstyle{plain}

\theoremstyle{plain}
\newtheorem{lemma}{Lemma}

\theoremstyle{remark}
\newtheorem*{rmk*}{Remark}

\theoremstyle{remark}
\newtheorem*{note*}{Note}

\def\pa{PA}
\def\z{\zeta}
\def\Z{\mathbb{Z}}
\def\R{\mathbb{R}}
\def\C{\mathbb{C}}
\def\cal{\mathcal}
\def\ds{\displaystyle}

\newcommand{\MP}[1]{\marginpar{\tiny #1}}
\newcommand{\parag}[1]{\medbreak\par\noindent{\bf\em#1}}
\newcommand{\Img}[2]{\includegraphics[width=#1truecm]{#2}}

\def\PSAW{prudent walk}
\def\PSAP{prudent polygon}
\def\W{W} 
\def\XX{X} 
\def\YY{Y} 
\def\ZZ{Z} 

\title{The 
Enumeration of Prudent Polygons by Area
and its Unusual Asymptotics}
\author{Nicholas R. Beaton, Philippe Flajolet, and Anthony J. Guttmann}

\date{November 29, 2010}

\begin{document}

\begin{abstract}
Prudent walks  are special self-avoiding  walks that never take a step
towards an already occupied  site, and \emph{$k$-sided  prudent walks}
(with $k=1,2,3,4$) are,  in essence, only   allowed to grow along  $k$
directions.  Prudent polygons are prudent walks that return to a point
adjacent to their starting point.    Prudent walks and polygons   have
been previously enumerated  by length and perimeter (Bousquet-M\'elou,
Schwerdtfeger; 2010).    We consider the enumeration  of \emph{prudent
polygons}  by \emph{area}.  For the  3-sided variety, we find that the
generating function is    expressed in terms  of a  $q$-hypergeometric
function, with  an    accumulation of  poles towards    the   dominant
singularity.  This expression reveals  an unusual asymptotic structure
of the number of polygons of area~$n$,  where the critical exponent is
the  transcendental number $\log_23$   and and the amplitude  involves
tiny oscillations.  Based on  numerical  data, we also expect  similar
phenomena to occur  for 4-sided polygons.   The asymptotic methodology
involves an original  combination of  Mellin transform techniques  and
singularity  analysis, which is  of potential interest  in a number of
other asymptotic enumeration problems.
\end{abstract}

\maketitle

\section{\bf Introduction}\label{sec:intro}

\def\SAP{\operatorname{SAP}}
\def\SAW{\operatorname{SAW}}

The problem of enumerating \emph{self-avoiding walks} (SAWs) and
\emph{polygons}
(SAPs) on a  lattice is a famous  one, 
whose complete solution has thus far   remained most elusive. 
For the square lattice, it is conjectured that the number
${\SAW}_n$ of walks  of \emph{length}~$n$ and the number 
${\SAP}_n$ of polygons of \emph{perimeter}~$n$ each satisfy an asymptotic formula of
the general form
\begin{equation}\label{eq0}
C\cdot \mu^n\cdot n^{\beta},
\end{equation}
where  $C,\mu\in\R_{>0}$ and $\beta\in\R$.  (In the case of polygons, it
is understood that~$n$ must be restricted to even values.)
The number~$\mu$ is the ``growth constant'' and the number~$\beta$ is often
referred to as the ``critical exponent''. 
More precisely, the following expansions are conjectured,
\begin{equation}\label{eq1}
{\SAW}_n \mathop{\sim}_{n\to\infty} C_1 \cdot
\mu^n \!\!{} \cdot   n^{11/32},
\qquad
{\SAP}_n \mathop{\sim}_{n\to\infty} C_2 \cdot 
\mu^n \!\!{} \cdot   n^{-5/2},
\end{equation}
for some $C_1,C_2>0$. 

For the square lattice, numerical methods based on acceleration of convergence and differential approximants suggest the value $\mu = 2.6381585303\ldots.$ This estimate is indistinguishable from the solution of the biquadratic equation $13\mu^4-7\mu^2-581 = 0$, which we consider to be a useful mnemonic. This was observed by Conway, Enting and Guttmann \cite{CoEnGu93} in 1993, and verified to 11 significant digits by Jensen and Guttmann \cite{JeGu99} in 2000, based on extensive numerical analysis of the sequence $({\SAP}_{2n})$ up to $2n=90$.
(Remarkably enough, for the
honeycomb lattice,   it had  long  been  conjectured  that the  growth
constant   of walks is the  biquadratic  number $\sqrt{2+\sqrt{2}}$, a
fact  rigorously established   only   recently  by  Duminil-Copin  and
Smirnov~\cite{DuSm10}.)

As    regards     critical    exponents,    the  conjectured     value
$\beta=\frac{11}{32}$ for  walks  is supported  by results  of Lawler,
Schramm,  and  Werner   that relate the    self-avoiding walk  to  the
``stochastic Loewner Evolution''  (SLE)  process of index $8/3$;  see,
for      instance,  the  account      in  Werner's   inspiring lecture
notes~\cite{Werner04}.    For     (unrooted)  polygons,    the   value
$\beta=-\frac52$  was  suggested  by numerical  analysis  of the exact
counting sequence,  with    an  agreement   to the   seventh   decimal
place~\cite{JeGu99}. It is also supported by the observation that many
simplified,   exactly    solvable,   (naturally rooted)    models   of
self-avoiding  polygons  appear   to exhibit  an  $n^{-3/2}$ universal
behaviour -- for these  aspects,      we  refer to      the survey  by
Bousquet-M\'elou    and    Brak~\cite{BoBr09},   as   well     as  the
books~\cite{Finch03,FlSe09,Rensburg00}.

Regarding  lattice polygons,  which are  closed  walks, there  is also
interest   in  enumeration  according   to \emph{area},  rather   than
perimeter. This question has  analogies  with the  classical  unsolved
problem of  enumerating polyominoes, also known  as animals  (these may
have  ``holes''),   according to the    number  of cells  they contain.
Conjecturally~\cite{JeGu00}, the number $a_n$ (respectively, $b_n$) of
polygons  (respectively, polyominoes)  comprised of~$n$  cells satisfy
asymptotic estimates of the form
\[
a_n\mathop{\sim}_{n\to\infty} C_3 \cdot (3.9709\ldots)^n\cdot  n^{-1},
\qquad
b_n\mathop{\sim}_{n\to\infty} C_4\cdot (4.0625\ldots)^n\cdot  n^{-1}, 
\]
for some $C_3,C_4\in\R_{>0}$; these asymptotic  estimates are still of
the   form~\eqref{eq0},  with   the  critical   exponent   $\beta=-1$.
Interestingly enough,  the  critical exponent~$\beta=0$, corresponding
to a simple pole of  the associated generating function, is  otherwise
known to arise  in  several simplified models, such  as column-convex, \MP{True?}
convex, and directed polyominoes~\cite{Bousquet96b,BoBr09,FlSe09}.

As the foregoing discussion  suggests, there is  considerable interest
in solving, exactly,  probabilistically, or asymptotically, restricted
models of self-avoiding walks and polygons.  Beyond serving to develop
informed  conjectures regarding more complex   models, this  is  relevant to
areas  such  as statistical  physics and  the statistical mechanics of
polymers~\cite{Rensburg00}. For combinatorialists, we may observe that
consideration  of such  models has served   as a powerful incentive to
develop       new  counting       methods   based  on       generating
functions~\cite{BoBr09,DeVi84,FlSe09,Temperley81},  including transfer
matrix methods and what is known as the ``kernel method''.

The  present   article focuses on   a   special type  of self-avoiding
polygons, the \emph{3-sided prudent  polygons} (to be defined shortly -- see
Definition~\ref{prud-def} in Section~\ref{prud-sec}),
when  these are enumerated according to  area.  
Roughly, a walk is prudent if it never takes a step towards an already 
occupied site and it is 2-, 3-, 4-sided if it has,
respectively, 2, 3, or 4 allowed directions of growth;
a prudent polygon is a prudent walk
that is almost closed.
For  area~$n$,  we will obtain a
precise asymptotic formula (Theorem~\ref{thm-asy0} below),
\begin{equation}\label{eq2}
\pa_n\sim 
C(n)\cdot 2^n\cdot n^g,
\end{equation}
one that    has several  distinguishing  features:  $(i)$~the critical
exponent is the   \emph{transcendental number} $g=\log_23$,   in sharp
contrast with   previously  known examples  where  it is  invariably a
``small'' rational number; $(ii)$~the   multiplier~$C$ is no  longer a
constant, but a bounded quantity    that \emph{oscillates} around the   value
$0.10838\ldots$ and does so with  a minute amplitude of $10^{-9}$. The
oscillations cannot be revealed by  any standard numerical analysis of
the  counting  sequence~$\pa_n$,  but  such a   phenomenon may well be
present in other models, and, if so, it could change our whole view of
the asymptotic behaviour of such models.

\parag{Plan and results of the paper.}
Prudent walks are defined   in Section~\ref{prud-sec}, where   we also
introduce   the  2-sided,    3-sided,  and  4-sided  varieties.    The
enumeration of 2- and 3-sided  walks and polygons by  \emph{perimeter}
is      the      subject    of          insightful      papers      by
Bousquet-M\'elou~\cite{Bousquet10}                                 and
Schwerdtfeger~\cite{Schwerdtfeger10} who obtained     both       exact
generating function expressions  and precise  asymptotic results.   In
Subsections~\ref{2prud-subsec}--\ref{4prud-subsec},
 we   provide the algebraic derivation of the corresponding area results: the enumeration of \emph{2- and 3- sided polygons according to  area} is treated there;
see Theorem~\ref{thm:3sided_gen_func} for  our first main result.  For
completeness, we also derive a functional  equation for the generating
function  of  4-sided prudent  polygons  (according   to area),  which
parallels                an             incremental       construction
of~\cite[\S6.5]{Bousquet10} -- this   functional equation suffices   to
determine    the   counting     sequence   in      polynomial    time.
Section~\ref{sec:asymptotics}, dedicated  to the asymptotc analysis of
the number of 3-sided prudent polygons, constitutes what  we feel to be
the   main   contribution    of   the  paper.     We  start     from a
$q$-hypergeometric representation   of   the generating   function  of
interest, $\pa(z)$, and  proceed to analyse  its singular structure: it
is found  that  $\pa(z)$  has  poles at  a  sequence  of  points  that
accumulate geometrically fast to~$\frac12$; then, the Mellin transform
technology~\cite{FlGoDu95} provides access to the asymptotic behaviour
of  $\pa(z)$, as $z\to\frac12$  in   extended regions of the   complex
plane.  Singularity analysis~\cite[Ch.VI--VII]{FlSe09} finally enables
us to determine the  asymptotic form of the coefficients~$\pa_n$
(see  Theorem~\ref{thm-asy0}) and   even derive  a complete asymptotic
expansion  (Theorem~\ref{thm-asy1}).   As already  mentioned,      the
non-standard character  of   the  asymptotic  phenomena  found  is   a
distinctive  feature.  Section~\ref{conclu-sec} concludes
the  paper with brief remarks  relating  to asymptotic methodology. In
particular, experiments suggest  that similar asymptotic phenomena are
likely to  be    encountered in the  enumeration   of  4-sided prudent
polygons.

A preliminary announcement of the results of the present paper is the object
of the communication~\cite{BeFlGu10}.

\section{\bf Prudent walks and polygons}\label{prud-sec}

One interesting   sub-class  of 
self-avoiding walks (SAWs)  for  which   a  number  of exact
solutions have been recently  found are \emph{prudent}  walks. 
Introduced by Pr\'{e}a \cite{Prea97},  these are SAWs which never take
a step towards an already occupied node. Exact solutions of {\PSAW}s on a
2-dimensional square lattice were later studied by Duchi
\cite{Duchi05}             and            Bousquet-M\'{e}lou
\cite{Bousquet10}, who    were  able   to  enumerate    certain
sub-classes. The enumeration of the corresponding class of polygons is
due  to Schwerdtfeger~\cite{Schwerdtfeger10}.   In   this section,  we
first  recall  the  definition of  prudent   walks and  polygons, then
summarize    the known  results   of~\cite{Bousquet10,Schwerdtfeger10}
relative to   their  enumeration   by   length or   perimeter;     see
Subsection~\ref{prud-subsec}, where 2-,  3-, and 4-sided prudent walks
are introduced.  We then examine the enumeration of polygons according
to area, in each of the three non-trivial  cases.  The case of 2-sided
polygons   is easy enough (Subsection~\ref{2prud-subsec}).    
The main result of
this   section  is Theorem~\ref{thm:3sided_gen_func}  of Subsection~\ref{3prud-subsec},
which   provides  an explicit    generating    function for    3-sided
polygons -- it is  on  this expression  that  our subsequent asymptotic
treatment    is  entirely  based.  In     the case of   4-sided (i.e.,
``general'') polygons,   we  derive 
in Subsection~\ref{4prud-subsec} a system    of functional
equations  that determines the generating   function and amounts to  a
polynomial-time  algorithm   for the generation   of the counting
sequence.

\begin{figure}
 \begin{center}
	\subfloat[]{\includegraphics[scale=0.8]{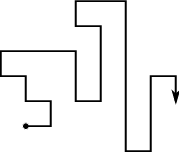}}
	\hspace{0.5cm}
	\subfloat[]{\includegraphics[scale=0.8]{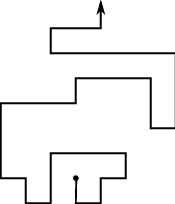}}
	\hspace{0.5cm}
	\subfloat[]{\includegraphics[scale=0.8]{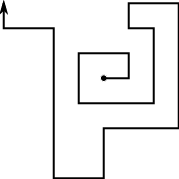}}
	\hspace{0.5cm}
	\subfloat[]{\includegraphics[scale=0.8]{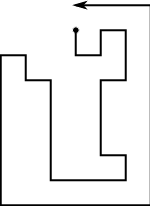}}
  \caption{\label{walkspoly-fig} Examples of (A) a two-sided prudent SAW; (B) a three-sided prudent SAW; (C) an (unrestricted) prudent SAW; (D) a prudent SAW leading to a prudent SAP.}
  \end{center}
\end{figure}

\subsection{Main definitions and results.}\label{prud-subsec}
We use  the same classification scheme  as the authors
of~\cite{Bousquet10,Schwerdtfeger10}. By definition, 
the endpoint of every {\PSAW} always lies on the boundary of the smallest
lattice rectangle which contains the entire  walk, referred to here as
the \emph{bounding  box} or just \emph{box}. This  property leads to a
natural classification of {\PSAW}s (see Figure~\ref{walkspoly-fig}).

\begin{definition}\label{prud-def}
Let $\omega$ be a prudent walk of length $n$, and let $\omega_i$ be the prudent walk comprising the first $i$ steps of $\omega$. Let $b_i$ be the bounding box of $\omega_i$. Then $\omega$ is \emph{1-sided} if $\omega_i$ ends on the north side of $b_i$ for each $i=0,1,\ldots,n$; \emph{2-sided} if each $\omega_i$ ends on the north or east sides of $b_i$; \emph{3-sided} if each $\omega_i$ ends on the north, east or west sides of $b_i$ (with one caveat, described below); and \emph{4-sided} (or \emph{unrestricted}) if each $\omega_i$ may end on any side of $b_i$.
\end{definition}

\begin{figure}
 \begin{center}
	\subfloat[]{\includegraphics[scale=0.8]{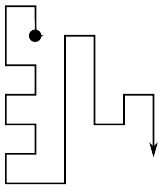}}
	\hspace{2cm}
	\subfloat[]{\includegraphics[scale=0.8]{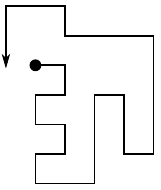}}
  \caption{\label{pathological-fig} Examples of the (A) prudent walks and (B) prudent polygons which we exclude from the definition of 3-sided.}
  \end{center}
\end{figure}

\begin{rmk*}
The  issue  with  3-sided  {\PSAW}s  is encapsulated  by  the walk $(0,0)
\rightarrow (1,0) \rightarrow (1,-1) \rightarrow (0,-1)$. If one draws
the walk's box  after each (discrete) step,  then it is clear that the
walk  always  ends  on the   north,  east  or  west sides,   seemingly
fulfilling  the 3-sided requirement. However, if  the walk is taken to
be continuous, then along  the step  $(1,-1) \rightarrow (0,-1)$,  the
endpoint is only on the south side. In general this occurs when a walk
steps from the south-east corner  of its box  to the south-west corner
(or vice versa) when the box has width one and non-zero height.  Allowing such walks forces us to account for structures like those in Figure~\ref{pathological-fig}; while this is certainly possible, it complicates a number of rational terms and contributes little to the asymptotic behaviour of the model. For this reason we follow the examples of Bousquet-M\'{e}lou and Schwerdtfeger \cite{Bousquet10, Schwerdtfeger10} and exclude these cases.

An equivalent definition of the 3-sided walks considered here is to forbid two types of steps: when the box has non-zero width, a south step may not be followed by a west (resp. east) step when the walk is on the east (resp. west) side of the box.
\end{rmk*}

\parag{Walks (length).}
The  enumeration  of 1-sided {\PSAW}s   (also known as \emph{partially
directed walks}) is  straightforward, and the  generating function for
such walks is rational -- we will not discuss  these any further. Duchi
\cite{Duchi05} successfully  found the generating function for 2-sided
walks,    showing    it   to      be    algebraic.  Bousquet-M\'{e}lou
\cite{Bousquet10} solved the problem of  3-sided {\PSAW}s, finding the
generating  function to be  non-D-finite.  For the unrestricted  case,
functional equations were found by  both Duchi and Bousquet-M\'{e}lou,
but at present these equations remain unsolved.

The  dominant singularity of  the generating functions for 2-sided and
3-sided  {\PSAW}s  is  a  simple pole at   $\rho  = 0.4030317168...$, the
smallest   real  root   of   $1-2x-2x^2+2x^3$.  Dethridge    and  Guttmann
\cite{DeGu08} conjecture that the same is true for unrestricted {\PSAW}s,
based on a computer-generated series  of 100 terms. They  also conjectured
that  the    generating   function     for  unrestricted     {\PSAW}s  is
non-holonomic\footnote{ A function  is said to  be \emph{holonomic} or
\emph{D-finite} if   it   is the solution  to   a  linear differential
equation   with      polynomial   coefficients;          see
\cite[Ch.~6]{Stanley99} and~\cite[Ap.~B.4]{FlSe09}.}  (non-D-finite). 
Here is a summary of known results, with the estimates tagged with a question mark~(?)
being conjectural ones:
\begin{equation}\label{walk-tab}
\hbox{\small
\begin{tabular}{lcll}
\hline\hline 
 \em Walks & \em Generating function & 
\multicolumn{1}{c}{\em Asymptotic number} & \em References
\\
\hline
2-sided & algebraic & $\ds \kappa_2\cdot \rho^{-n},\quad \rho^{-1}\simeq
2.481$ & Duchi~\cite{Duchi05},
Bousquet-M\'elou~\cite{Bousquet10}
\\
3-sided & non-holonomic & $\ds \kappa_3\cdot \rho^{-n},\quad 
	\rho^{-1}\simeq2.481$ & Bousquet-M\'elou~\cite{Bousquet10}
\\
4-sided & functional equation & $\ds \kappa_4\cdot \rho^{-n},\quad 
\rho^{-1}\simeq2.481$ (?) & Dethridge \& Guttmann~\cite{DeGu08}.
\\
\hline\hline
\end{tabular}}
\end{equation}
The     values   of    the     multipliers,  after~\cite{DeGu08},  are
$\kappa_2=2.51...$,        $\kappa_3=6.33...$      and     (estimated)
$\kappa_4\approx 16.12$.

\parag{Prudent polygons (perimeter).}
Self-avoiding  polygons (SAPs) are self-avoiding walks  which end at a
node adjacent to  their starting  point (excluding  walks of a  single
step). If the walk has length $n-1$ then the polygon is said to have \emph{perimeter} $n$.

\begin{definition}\label{prud-poly-def}
A \emph{prudent} self-avoiding polygon ({\PSAP}) is a SAP for
which the  underlying SAW  is prudent.   In the  same way, a   {\PSAP} is
1-sided  (resp. 2-sided, etc.)   if its  underlying  {\PSAW}   is 1-sided
(resp. 2-sided, etc.).
\end{definition}

A 1-sided {\PSAP}   starting at    $(0,0)$  must end  at $(0,1)$,    thus
consisting only of   a single  row of  cells  and having  a   rational
generating  function. The   enumeration of 2-   and  3-sided {\PSAP}s  by
perimeter       has     been         addressed     by    Schwerdtfeger
\cite{Schwerdtfeger10}. The non-trivial  2-sided {\PSAP}s  are essentially
inverted   bargraphs  \cite{PrBr95}, and so  the   2-sided case has an
algebraic generating function.  Schwerdtfeger finds  the 3-sided {\PSAP}s
to have a non-D-finite generating function.

If  $PP^{(k)}(z) =   \sum_{n\geq0}  p^{(k)}_n  z^n$   is defined  to  be  the
half-perimeter generating  function for $k$-sided prudent polygons (so
$p^{(k)}_n$ is  the number of  $k$-sided  prudent polygons with perimeter
$2n$), then the following holds~\cite{Schwerdtfeger10}:
$(i)$~the   dominant singularity   of  $PP^{(2)}(z)$  is   a square   root
singularity  at   $\sigma = 0.2955977...$,  the   unique real  root of
$1-3x-x^2-x^3$. So $p^{(2)}_n \sim \lambda_2 \sigma^{-n} n^{-3/2}$ as $n\to\infty$,
where $\lambda_2$ is a constant;
$(ii)$~the dominant singularity of $PP^{(3)}(z)$ is a 
square root singularity at $\tau=0.24413127...$, where $\tau$ is 
the unique real root of 
$\tau^5+6\tau^3-4\tau^2+17\tau-4$.
So $p^{(3)}_n \sim  
\lambda_3\tau^{-n} n^{-3/2}$ as $n\to\infty$, where $\lambda_3$ is a constant.
Schwerdtfeger has furthermore classified 4-sided {\PSAP}s in such a way as to allow for functional equations in the generating functions to be written. Unfortunately no one
 has thus far been able to obtain a solution from said equations.
We again present a summary of known results.

\begin{equation}\label{poly-tab}
\hbox{\small
\begin{tabular}{lcll}
\hline\hline 
 \em Polygons (perimeter)& \em Generating function 
& \multicolumn{1}{c}{\em Asymptotic number} & \em References
\\
\hline
2-sided & algebraic & $\ds \lambda_2\cdot \sigma^{-n}
n^{-3/2},\quad \sigma^{-1}\simeq 3.382$ & Schwerdtfeger~\cite{Schwerdtfeger10}
\\
3-sided & non-holonomic & $\ds \lambda_3\cdot \tau^{-n}n^{-3/2},\quad \tau^{-1}\simeq4.096$ & Schwerdtfeger~\cite{Schwerdtfeger10}
\\
4-sided & functional equation & 
$\ds \lambda_4\cdot \upsilon^{-n} n^{-\delta},
\hbox{\hspace*{2.5truemm}}\quad \upsilon^{-1}\approx4.415$ (?)  
 & Dethridge \emph{et al.}~\cite{DeGaGuJe09}
\\
\hline\hline
\end{tabular}
}
\end{equation}
The empirical estimates regarding 4-prudent polygons are 
taken from~\cite{DeGaGuJe09}. They are somewhat imprecise,
and it is suspected that \MP{check?}
the critical exponent satisfies $\delta=-3.5\pm0.1$,
with $\delta=-7/2$ a compatible value.

\parag{Prudent polygons (area).}
The focus of this paper  is on the enumeration of  {\PSAP}s by \emph{area}, rather
than perimeter. The constructions we use here are essentially the same
as  Schwerdtfeger's \cite{Schwerdtfeger10}; the  resulting  functional  equations  and  their solutions, however, turn out to be quite different,  as will be revealed by
the peculiar singularity structure of the generating functions and the
non-trival asymptotic form of the coefficients. We  have  modified Schwerdtfeger's
construction for  4-sided   {\PSAP}s  slightly to  allow   for an  easier
conversion into a recursive form (see Subsection \ref{4prud-subsec}).

We will denote the area  generating function for $k$-sided {\PSAP}s by
$\pa^{(k)}(q) = \sum_{n\geq1} \pa^{(k)}_n q^n$. For 3-  and 4-sided {\PSAP}s, it
is necessary  to measure more  than just the area   -- in these cases,
additional    \emph{catalytic}    variables    will   be  used    (see
\cite{Bousquet10} for a more thorough explanation).

\subsection{Enumeration of 2-sided polygons by area.}\label{2prud-subsec}
The non-trivial 2-sided {\PSAP}s can be constructed from bargraphs. 
Let $B(q) = \sum_{n\geq1} b_n q^n$ be the area generating function for these objects.
The area generating function for bargraphs, $B(q)$, is
\[B(q) = \frac{q}{1-2q}\]
and so $b_n = 2^{n-1}$ for $n\geq 1$. (Bargraphs are a graphical representation
of integer compositions.)
%

%



\begin{prop}\label{prop:2sided_gf}
The area generating function for 2-sided {\PSAP}s is
\[\pa^{(2)}(q) = \frac{2q}{1-2q} + \frac{2q}{1-q},\]
and so the number of such polygons is $\pa^{(2)}_n = 2^{n}+2$ for $n\geq 1$.
\end{prop}

\begin{proof}
A 2-sided prudent polygon must end at either $(0,1)$ or $(1,0)$. Reflection in the line $y=x$ will not invalidate the 2-sided property, so it is sufficient to enumerate those polygons ending at $(1,0)$ and then multiply the result by two.

The underlying 2-sided {\PSAW} cannot step above the line $y=1$, nor to any point $(x,y)$ where $x,y < 0$. So any polygon beginning with a west step must be a single row of cells to the left of the $y$-axis. The generating functions for these polygons is then $\frac{q}{1-q}$.

A polygon starting with a south or east step must remain on the east side of its box until it reaches the line $y=1$, at which point it has no choice but to take west steps back to the $y$-axis. It can hence be viewed as an upside-down bargraph with north-west corner $(0,1)$. The area generating function for these objects is $B(q) = \frac{q}{1-2q}$.

Adding these two possibilities together and doubling gives the result.
\end{proof}

\subsection{Enumeration of 3-sided polygons by area.}\label{3prud-subsec}
When constructing 3-sided {\PSAP}s, we will use a single catalytic variable which measures width. To do so we will need to measure bargraphs by width. Let $B(q,u) = \sum_{n\geq1}\sum_{i\geq1}b_{n,i}q^n u^i$ be the area-width generating function for bargraphs (so $b_{n,i}$ is the number of bargraphs with area $n$ and width $i$).

The area-width generating function for bargraphs, $B(q,u)$, satisfies the equation
\begin{equation}\label{eqn:bargraph_func_eq_width}
B(q,u) = \frac{qu}{1-q} + \frac{qu}{1-q}B(q,u),
\end{equation}
which is obtained by successively adding columns. Accordingly,
by solving the functional equation, we obtain
\[B(q,u) = \frac{qu}{1-q-qu}\]
and so $b_{n,i} = \binom{n-1}{i-1}$ for $n,i \geq 1$.
(Clearly, $b_{n,i}$ counts compositions of~$n$ into~$i$ summands.)


%
%

Let ${\W}(q,u) = \sum_{n\geq1}\sum_{i\geq1}w_{n,i}q^n u^i$ be the area-width generating function for 3-sided {\PSAP}s which end at $(-1,0)$ in a counter-clockwise direction. As we will see, this is the most complex type of 3-sided {\PSAP}; everything else is either a reflection of this or can be constructed from something simpler.

\begin{figure}
\centering
\subfloat{\includegraphics[scale=1.5]{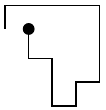}}
\hspace{0.5cm}
\subfloat{\includegraphics[scale=0.75]{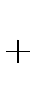}}
\hspace{0.5cm}
\subfloat{\includegraphics[scale=1.5]{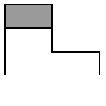}}
\hspace{0.5cm}
\subfloat{\includegraphics[scale=0.75]{img_plus_sign}}
\hspace{0.5cm}
\subfloat{\includegraphics[scale=1.5]{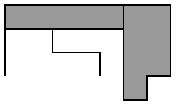}}
\caption{The decomposition used to construct 3-sided prudent polygons.}
\label{fig:3sided_construction_decomp}
\end{figure}

\begin{lemma}\label{lem:-10_func_eq_width}
The area-width generating function for 3-sided prudent polygons ending at $(-1,0)$ in a counter-clockwise direction, ${\W}(q,u)$, satisfies the functional equation 
\begin{equation}\label{eqn:-10_func_eq_width}
{\W}(q,u) = qu(1+B(q,u)) + \frac{q}{1-q}({\W}(q,u)-{\W}(q,qu)) + qu(1+B(q,u)){\W}(q,qu)
\end{equation}
\end{lemma}

\begin{proof}
The underlying {\PSAW} cannot step to any point $(x,y)$ with $x,y<0$, nor to any point with $x<-1$. It must approach the final node $(-1,0)$ from above. So the only time the endpoint can be on the west side of the box and not the north or south is when the walk is stepping south along the line $x=-1$. So prior to reaching the line $x=-1$, the walk must in fact be 2-sided. Note that the north-west corner of the box must be a part of of the polygon.

If the walk stays on or below the line $y=1$, then (as will be seen in Proposition \ref{prop:2sided_gf}), it either reaches the point $(0,1)$ with a single north step, or by forming an upside-down bargraph. This must then be followed by a west step to $(-1,1)$, then a south step. This will form either a single square or a bargraph with a single square attached to the north-west corner, giving the first term on the right-hand side of (\ref{eqn:-10_func_eq_width}).

Since the north-west corner of the box of any of these polygons is part of the polygon, it is valid to add a row of cells to the top of an existing polygon (so that the west sides line up). This can be done to any polygon. If the new row is not longer than the width of the existing polygon we obtain the term
\[\sum_{n\geq 1}\sum_{i\geq1}w_{n,i}q^n u^i \cdot \sum_{k=1}^i q^k = q \sum_{n\geq 1}\sum_{i\geq1}w_{n,i}q^n u^i \cdot\frac{1-q^i}{1-q},\]
giving the second term in the right-hand side of (\ref{eqn:-10_func_eq_width}).
\begin{note*} For the remainder of this subsection, we will omit unwieldy double or triple sums like the one above, and instead give recursive relations only in terms of the generating functions.
\end{note*}

Instead, the new row may be longer than the width of the existing polygon. In this case, as the walk steps east along this new row, it will reach a point at which there are no occupied nodes south of its position, and it will hence be able to step south in a prudent fashion. It must then remain on the east side of the box until reaching the north side, at which point it steps west to $x=-1$ and then south to the endpoint. This effectively means we have added a row of length equal to the width $+1$, and then (possibly) an arbitrary bargraph. So we obtain
\[qu{\W}(q,qu)(1+B(q,u))\]
which gives the final term in the right-hand side of (\ref{eqn:-10_func_eq_width}).
\end{proof}

\begin{lemma}\label{lem:-10_gen_func}
The area-width generating function for 3-sided {\PSAP}s ending at $(-1,0)$ in a counter-clockwise direction is
\[{\W}(q,u) = \sum_{m=0}^{\infty}F(q,q^mu)\prod_{k=0}^{m-1}G(q,q^ku),\]
where
\[
F(q,u)  = \frac{qu(1-q)^2}{(1-2q)(1-q-qu)},\qquad 
G(q,u)  = \frac{-q(1-q-u+qu-q^2u)}{(1-2q)(1-q-qu)}.
\]
\end{lemma}

\begin{proof}
Substituting $BW(q,u)=\frac{qu}{1-q-qu}$ into (\ref{eqn:-10_func_eq_width}) and rearranging gives
\begin{equation}\label{eqn:-10_func_eq_FG}
{\W}(q,u) = F(q,u)+G(q,u){\W}(q,qu)
\end{equation}
Substituting $u \mapsto uq$ gives
\begin{equation}\label{eqn:-10_func_eq_FG_q2}
{\W}(q,qu)=F(q,qu)+G(q,qu){\W}(q,q^2u)
\end{equation}
and combining these yields
\begin{equation}\label{eqn:-10_func_eq_FG_q2sub}
{\W}(q,u)=F(q,u)+F(q,qu)G(q,u)+G(q,u)G(q,qu){\W}(q,q^2u).
\end{equation}
Repeating for $u \mapsto q^2u, q^3u,..., q^Mu$ will give
\begin{equation}\label{eqn:-10_func_eq_FG_gMsub}
{\W}(q,u) = \sum_{m=0}^{M}F(q,q^mu)\prod_{k=0}^{m-1}G(q,q^ku)+\prod_{m=0}^{M}G(q,q^mu){\W}(q,q^{M+1}u).
\end{equation}

We now seek to take $M\rightarrow \infty$. To obtain the result stated in the Lemma, it is necessary to show that 
\[\sum_{m=0}^{M}F(q,q^mu)\prod_{k=0}^{m-1}G(q,q^ku)\]
converges, and 
\[\prod_{m=0}^{M}G(q,q^mu){\W}(q,q^{M+1}u) \rightarrow 0\]
as $M\to\infty$ (both considered as power series in $q$ and $u$).

Both $F$ and $G$ are 
bivariate power series in $q$ and $u$.
We have that
\begin{align*}
F(q,u) & = qu + q^2(u+u^2)+q^3(2u+2u^2+u^3)+O(q^4)\\
G(q,u) & = q(-1+u)+q^2(-2+u+u^2)+q^3(-4+2u+2u^2+u^3) + O(q^4)
\end{align*}
It follows that $F(q,q^mu) = O(q^{m+1})$ and $G(q,q^ku) = O(q)$ for all $m,k \geq 0$. So we have that
\[F(q,q^mu)\prod_{k=0}^{m-1}G(q,q^ku) = O(q^{2m+1})\]
So considered as a power series in $q$ and $u$, the first term in the right-hand side of (\ref{eqn:-10_func_eq_FG_gMsub}) does converge to a fixed power series as $M\to\infty$.

By the same argument, we have that 
\[\prod_{m=0}^{M}G(q,q^mu)\to 0\]
as $M\to\infty$. So it suffices to show that ${\W}(q,q^{M+1}u)$ converges to a fixed power series. But now every term in the series ${\W}(q,u)$ has at least one factor of $u$ (since every polygon has positive width), so it immediately follows that ${\W}(q,q^{M+1}u)\to0$ as $M\to\infty$.

So both terms in (\ref{eqn:-10_func_eq_FG_gMsub}) behave as required as $M\to\infty$, and the result follows.
\end{proof}

\begin{theorem}\label{thm:3sided_gen_func}
The area generating function for 3-sided prudent polygons is
\begin{align*}\begin{split}
\pa^{(3)}(q) & = \frac{-2q^3(1-q)^2}{(1-2q)^2}\sum_{m=1}^{\infty}\frac{(-1)^m q^{2m}}{(1-2q)^m(1-q-q^{m+1})} \prod_{k=1}^{m-1}\frac{1-q-q^k+q^{k+1}-q^{k+2}}{1-q-q^{k+1}}\\&\,\hspace{9.5cm}+\frac{2q(3-10q+9q^2-q^3)}{(1-2q)^2(1-q)}\\
& = 6q+10q^2+20q^3+42q^4+92q^5+204q^6+454q^7+1010q^8+2242q^9+4962q^{10}...\end{split}
\end{align*}
\end{theorem}

\begin{proof}
A 3-sided {\PSAP} must end at $(-1,0), (0,1)$ or $(1,0)$, in either a clockwise or counter-clockwise direction. Setting $u=1$ in ${\W}(q,u)$ gives the area generating function
\begin{equation}\label{eqn:-10_gen_func_u=1}
{\W}(q,1) = \frac{-q^3(1-q)^2}{(1-2q)^2}\sum_{m=1}^{\infty}\frac{(-1)^m q^{2m}}{(1-2q)^m(1-q-q^{m+1})}\prod_{k=1}^{m-1}\frac{1-q-q^k+q^{k+1}-q^{k+2}}{1-q-q^{k+1}} + \frac{q(1-q)^2}{(1-2q)^2}.
\end{equation}

A clockwise polygon ending at $(-1,0)$ can only be a single column, which has generating function 
\begin{equation}\label{eqn:3sided_single_col}
\frac{q}{1-q}.
\end{equation}

A counter-clockwise polygon ending at $(0,1)$ cannot step left of the $y$-axis or above the line $y=1$. While it is below this line, it must remain on the east side of its box, and upon reaching the line $y=1$, it must step west to the $y$-axis. It must therefore be a bargraph, with generating function 
\begin{equation}\label{eqn:3sided_bargraph_part}
\frac{q}{1-2q}.
\end{equation}

A reflection in the $y$-axis converts a polygon ending at $(-1,0)$ to one ending at $(1,0)$ in the opposite direction, and reverses the direction of a polygon ending at $(0,1)$. So adding together and doubling (\ref{eqn:-10_gen_func_u=1}), (\ref{eqn:3sided_single_col}) and (\ref{eqn:3sided_bargraph_part}) will cover all possibilities, and gives the stated result.
\end{proof}

\subsection{Enumeration of 4-sided polygons by area.}\label{4prud-subsec}
This case is included for completeness, as the results
are not needed in our subsequent asymptotic analysis.
A 4-sided  {\PSAP} may end  at any of  $(0,1), (1,0),  (0,-1), (-1,0)$ in
either  a clockwise  or  counter-clockwise  direction. Reflection  and
rotation leads to  an 8-fold symmetry, so  it  suffices to count  only
those ending at $(-1,0)$  in a counter-clockwise direction.  We modify
Schwerdtfeger's sub-classification slightly.

Let ${\XX}(q,u,v)  = \sum_{n\geq1}\sum_{i\geq1}\sum_{j\geq1}x_{n,i,j}q^n
u^i   v^j$  be the   generating function   for  those  polygons (class
$\mathcal{X}$)  for which  removing the  top  row does not change  the
width or  leave  two or more  disconnected  pieces, with $q$ measuring
area, $u$ measuring width and $v$ measuring height.

Let  ${\YY}(q,u,v) = \sum_{n\geq1}\sum_{i\geq1}\sum_{j\geq1}y_{n,i,j}q^n
u^i  v^j$ be the  generating function for  the  unit square plus those
polygons (class $\mathcal{Y}$) not in $\mathcal{X}$ for which removing
the rightmost column does not  change the height or  leave two or more
disconnected pieces, with $q$ measuring area, $u$ measuring height and
$v$ measuring width.

Let ${\ZZ}(q,u,v)  = \sum_{n\geq1}\sum_{i\geq1}\sum_{j\geq1}z_{n,i,j}q^n
u^i  v^j$   be  the  generating function   for   the  polygons  (class
$\mathcal{Z}$)  not  in   $\mathcal{X}$  or  $\mathcal{Y}$, with   $q$
measuring area, $u$ measuring width $-1$ and $v$ measuring height.

\begin{figure}
\centering
\renewcommand{\arraystretch}{0.3}
\begin{tabular}{l}
\raisebox{15pt}{$\cal X$~:} \qquad
\includegraphics[scale=1.25]{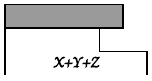}
\\
\raisebox{25pt}{$\cal Y$~:} \qquad
\subfloat{\includegraphics[scale=1.25]{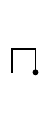}}
\hspace{0.3cm}
\subfloat{\includegraphics[scale=0.75]{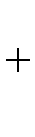}}
\hspace{0.3cm}
\subfloat{\includegraphics[scale=1.25]{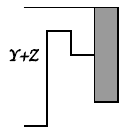}}
\hspace{0.3cm}
\subfloat{\includegraphics[scale=0.75]{img_plus_sign2}}
\hspace{0.3cm}
\subfloat{\includegraphics[scale=1.25]{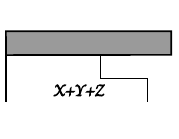}}
\\
\raisebox{30pt}{$\cal Z$~:} \qquad
\subfloat{\includegraphics[scale=1.25]{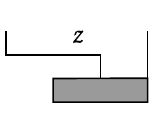}}
\hspace{0.3cm}
\raisebox{3pt}{\subfloat{\includegraphics[scale=0.75]{img_plus_sign2}}}
\hspace{0.3cm}
\subfloat{\includegraphics[scale=1.25]{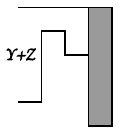}}
\end{tabular}

\caption{\label{XYZ-fig}
The decompositions used to construct 4-sided prudent
polygons in $\cal X, \cal Y, \cal Z$ (from top to bottom).}
\end{figure}

\begin{prop}\label{prop:4sided_func_eq}
The generating functions ${\XX}(q,u,v)$, ${\YY}(q,u,v)$ and ${\ZZ}(q,u,v)$ satisfy the functional equations
\begin{equation}\label{eqn:4sided_X_func_eq}
\begin{split}
{\XX}(q,u,v) = \frac{qv}{1-q}[{\XX}(q,u,v)-{\XX}(q,qu,v)]+\frac{qv}{1-q}[{\YY}(q,v,u)-{\YY}(q,v,qu)]\\+\frac{quv}{1-q}[{\ZZ}(q,u,v)-q{\ZZ}(q,qu,v)]
\end{split}
\end{equation}
\begin{equation}\label{eqn:4sided_Y_func_eq}
\begin{split}
{\YY}(q,u,v) = quv+\frac{qv}{1-q}[{\YY}(q,u,v)-{\YY}(q,qu,v)]+\frac{qv^2}{1-q}[{\ZZ}(q,v,u)-{\ZZ}(q,v,qu)]\\+quv[{\XX}(q,qv,u)+{\YY}(q,u,qv)+qv{\ZZ}(q,qv,u)]
\end{split}
\end{equation}
\begin{equation}\label{eqn:4sided_Z_func_eq}
{\ZZ}(q,u,v)=\frac{qv}{1-q}[{\ZZ}(q,u,v)-{\ZZ}(q,qu,v)]+qv{\YY}(q,qv,u)+quv{\ZZ}(q,u,qv)
\end{equation}
The generating function for 4-sided {\PSAP}s is then given by
\begin{align*}
PA^{(4)}(q) &= 8[X(q,1,1)+Y(q,1,1)+Z(q,1,1)]\\
&= 8q + 16q^2 + 40q^3 + 96q^4 + 232q^5 + 560q^6 + 1336q^7 + 3176q^8 + 7480q^{9} + 17528q^{10} \ldots
\end{align*}
\end{prop}

\begin{proof}
As with the 3-sided polygons in Lemma \ref{lem:-10_func_eq_width}, the walk cannot visit any point $(x,y)$ with $x,y<0$ or with $x<-1$. The walk must approach $(-1,0)$ from above, and must do so immediately upon reaching the line $x=-1$. So every polygon contains the north-west corner of its box. As in the 3-sided case, this leads to a construction involving adding rows to the top of existing polygons.

By definition, a polygon in $\mathcal{X}$ of width $i$ can be constructed by adding a row of length $\leq i$ to the top of any polygon of width $i$. Adding a row to a polygon in $\mathcal{X}$ gives 
\[\sum_{n\geq1}\sum_{i\geq1}\sum_{j\geq1}x_{n,i,j}q^n u^i v^j \cdot v\sum_{k=1}^{i}q^k = qv\sum_{n\geq1}\sum_{i\geq1}\sum_{j\geq1}x_{n,i,j}q^n u^i v^j \cdot\frac{1-q^i}{1-q},\]
which is the first term in the right-hand side of (\ref{eqn:4sided_X_func_eq}). Performing similar operations for polygons in $\mathcal{Y}$ and $\mathcal{Z}$ gives the rest of (\ref{eqn:4sided_X_func_eq}).
\begin{note*} Again, for the remainder of this subsection we give recursive relations purely in terms of the generating functions.
\end{note*}

Polygons not in $\mathcal{X}$ must also contain the north-east corner of their box. This leads to another construction involving adding columns to the right-hand side of existing polygons. To obtain a polygon in $\mathcal{Y}$ of height $i$, a new column of height $\leq i$ should be added to a polygon of height $i$ which contains the north-east corner of its box. So adding a column to a $\mathcal{Y}$ polygon gives
\[\frac{qv}{1-q}Y(q,u,v) - \frac{qv}{1-q}Y(q,qu,v)\]
which is the second term in the right-hand side of (\ref{eqn:4sided_Y_func_eq}). Performing a similar operation for $\mathcal{Z}$ polygons gives the third term in (\ref{eqn:4sided_Y_func_eq}).

Adding a new column to a polygon in $\mathcal{X}$ containing its north-east corner can be viewed as adding a sequence of rows on top of one another, and so if the new column has height $\geq2$ then the resulting polygon is actually in $\mathcal{X}$. If the new column has height one, however, the resulting polygon is in $\mathcal{Y}$. Isolating those polygons in $\mathcal{X}$ which contain their north-east corner is difficult; however, we can perform an equivalent construction by adding a row of length $i+1$ to any polygon of width $i$. Doing so to a polygon in $\mathcal{X}$ gives
\[quvX(q,qv,u)\]
and combining this with the same for $\mathcal{Y}$ and $\mathcal{Z}$ gives the fourth term in (\ref{eqn:4sided_Y_func_eq}). The $quv$ term is the unit square.

Polygons in $\mathcal{Z}$ also contain the south-east corner of their box. In a similar fashion to the constructions for $\mathcal{Y}$ and $\mathcal{Z}$, we can add a new row to the bottom of a polygon containing its south-east corner. To do so to a polygon in $\mathcal{Z}$ of width $i+1$ (remember $u$ measures width $-1$) requires a new row of width $\leq i$, so we obtain
\[\frac{qv}{1-q}Z(q,u,v)-\frac{qv}{1-q}Z(q,qu,v)\]
which is the first term in the right-hand side of (\ref{eqn:4sided_Z_func_eq}). 

Adding a new row to the bottom of something in $\mathcal{X}$ (containing its south-east corner) will give back something in $\mathcal{X}$, which will have been constructed by an alternate method described above. Adding a new row of length $\geq2$ to a polygon in $\mathcal{Y}$ will result in another polygon in $\mathcal{Y}$, which will also be constructible via alternate means. So we are left only with the possibility of adding a row of length one to the bottom of a polygon in $\mathcal{Y}$. This is analogous to the above description of adding a column of height one to the right of a polygon in $\mathcal{X}$; we now proceed by adding a column of height $i+1$ to a polygon in $\mathcal{Y}$ or $\mathcal{Z}$ of height $i$. Doing so gives the final two terms in (\ref{eqn:4sided_Z_func_eq}).
\end{proof}

\section{\bf Asymptotics}\label{sec:asymptotics}
For most lattice  object problems, finding  and solving the functional
equation(s) is the difficult part. Once a generating function has been
found,  the dominant  singularity is  often  quite obvious, and so the
asymptotic form  of  the coefficients  can  be easily described.   The
problem of  3-sided prudent polygons,  however, turns out to be rather
the  opposite.  The  functional equation (\ref{eqn:-10_func_eq_width})
was not terribly difficult to  obtain,  and its solution is  relatively
simple -- it only comprises a sum of products of rational functions of
$q$.

The  asymptotic  behaviour  of  this  model, on  the  other  hand,  is
considerably  more complex  than any  model we  have seen  before. The
dominant  singularity at $q=1/2$ is not even apparent from the representation
of Theorem~\ref{thm:3sided_gen_func}. As we shall see, there is in fact an accumulation of poles of the generating 
function\footnote{ Throughout this section
only dedicated to 3-sided prudent polygons, we omit redundant superscripts
and let $\pa_n$ and $\pa(q)$ represent, respectively, what was denoted
by $\pa^{(3)}_n$ and $\pa^{(3)}(q)$ 
in Section~\ref{prud-sec}.} $\pa(z)$ towards $q=1/2$. 
Accordingly, the nature of the dominant singularity at $q=1/2$ is rather
unusual: a singular expansion as $q$ approaches $1/2$ can be determined, but it involves
periodic fluctuations, a strong divergence from the standard simple type 
$Z^\alpha(\log Z)^\beta$, where $Z:=1-z/\rho$, with~$\rho$ (here equal
to $1/2)$ the dominant singularity of the generating function under consideration.
This is revealed by a Mellin analysis of $\pa(z)$  near its singularity, and the periodic fluctuations, which 
appear to be in a logarithmic scale, eventually echo  the geometric speed with which poles accumulate
at $1/2$. Then, thanks to a suitable extension to the complex plane,
the singular expansion can be transfered to coefficients
by the method known as singularity analysis~\cite[Ch.~VI]{FlSe09}. 
The net result is, for the
coefficients $\pa_n$, an asymptotic form that involves a standard element $2^n n^g$,
but multiplied by a \emph{periodic function} in~$\log_2 n$. The presence of oscillations, 
the transcendental character of the exponent~$g=\log_23$, and the minute amplitude of
these oscillations, about $10^{-9}$, are noteworthy features of this asymptotic problem.

\begin{theorem}\label{thm-asy0}
The number
 $\pa_n\equiv \pa^{(3)}_n$ 
of 3-sided prudent polygons of area~$n$ satisfies the estimate
\begin{equation}\label{pnasy}
\pa_n=
 \left[\kappa_0+\kappa(\log_2 n)\right]\,  2^n \cdot  n^{g}
+O\left(2^n\cdot n^{g-1}\log n\right),\qquad n\to\infty,
\end{equation}
where  the critical exponent is 
\[
g=\log_23\doteq1.58496\]
and the ``principal'' constant is
\begin{equation}\label{kap0}
\kappa_0= 
\frac{\pi}{9\log 2\, \sin(\pi g)\, \Gamma(g+1)}
\prod_{j=0}^\infty \frac{(1-\frac13 2^{-j})(1-\frac32 2^{-j})}{(1-\frac12 2^{-j})^2}
\doteq0.10838\,42946.
\end{equation}
The function $\kappa(u)$ is a smooth periodic
function of~$u$, with period~1, mean value zero, and 
amplitude 
$\doteq 1.54623\cdot 10^{-9}$, which is determined by its Fourier series representation:
\[
\kappa(u)=\sum_{k\in\Z\setminus\{0\}} \kappa_k e^{2ik\pi u}, \qquad 
\hbox{with}\quad \kappa_k=\kappa_0\cdot \frac{\sin (\pi g)}{\sin(\pi g+2ik\pi^2/\log2)}\cdot 
\frac{\Gamma (1+g)}{\Gamma(1+g+2ik\pi/\log2)}.
\]
\end{theorem}

The proof of the theorem occupies the next subsections, whose
organization reflects the informal description given above. We shall then discuss
the fine structure of subdominant terms in the asymptotic expansion of~$\pa_n$; cf Theorem~\ref{thm-asy1}.
Some quantities that appear repeatedly throughout this section are
tabulated in Figure~\ref{tab-fig} for convenience.

\begin{figure}\small
\begin{center}
\begin{tabular}{lcl}
\hline\hline
\multicolumn{1}{c}{\em Quantity} & \multicolumn{1}{c}{\em at $q=1/2$} & 
\multicolumn{1}{c}{\em reference}\\
\hline
$\ds u=\frac{q}{1-q}$ & 1 & Eq.~\eqref{acuv3}\\
$\ds v=\frac{1-q+q^2}{1-q}$ & $\ds\frac32$ &  Eq.~\eqref{acuv3}\\
$\ds a=\frac{q^2}{1-q+q^2}$ & $\ds\frac13$ & Eq.~\eqref{defa}\\
$\ds \gamma=\frac{\log v}{\log 1/q}$ & $\log_2(3/2)$ & Eq.~\eqref{defgam}\\
$\ds C(q)=\frac{2q(3-10q+9q^2-q^3)}{(1-q)(1-2q)^2}$ &
 $\ds {}\sim \frac{1}{4(1-2q)^2}$ & Eq.~\eqref{acuv} and~\eqref{eqn:C_expansion_1/2}\\
$\ds A(q)=\frac{2q(1-q)^2}{(1-2q)^2}$ & 
 $\ds {}\sim \frac{1}{4(1-2q)^2}$ & Eq.~\eqref{acuv} and~\eqref{eqn:A_expansion_1/2}\\
\hline\hline
\end{tabular}
\end{center}

\caption{\label{tab-fig}\small
A table of some of the recurring quantities of Section~\ref{sec:asymptotics},
their reduction at $q=1/2$ and the relevant equations in the text.}
\end{figure}

\subsection{Resummations}
We start with a minor reorganization of the formula provided by Theorem~\ref{thm:3sided_gen_func}:
completion of the finite products that appear there
leads to the equivalent \emph{$q$--hypergeometric} form
\begin{equation}\label{basic}
\pa(q) = C(q)+A(q) \cdot Q(1;q)\cdot \sum_{n=1}^{\infty}\frac{(-1)^n q^{2n}}{(1-2q)^n}\cdot \frac{1}{Q(q^n;q)}.
\end{equation}
Here and throughout this section, the notations are
\begin{equation}\label{acuv}
C(q) := \frac{2q(3-10q+9q^2-q^3)}{(1-q)(1-2q)^2}, \qquad
A(q):= \frac{2q(1-q)^2}{(1-2q)^2},
\end{equation}
and
\begin{equation}\label{acuv2}
Q(z;q):= 
Q\left(z;q;u(q),v(q)\right), 
\qquad\hbox{where}\quad 
Q(z;q;u,v)=\frac{(vz;q)_\infty}{(quz;q)_\infty},
\end{equation}
with 
\begin{equation}\label{acuv3}
u(q) = \frac{q}{1-q}, \qquad
v(q) = \frac{1-q+q^2}{1-q}.
\end{equation}
In the definition of~$Q$, the notation
$(x;q)_n$ represents the usual $q$-Pochhammer symbol:
\[
(x;q)_n=(1-x)(1-qx)\cdots (1-xq^{n-1}).
\]

\begin{lemma}\label{lem-cont0}
The function $\pa(q)$ is analytic in the open disc $|q|<\sqrt{2}-1$,
where it admits the convergent \emph{$q$-hypergeometric} representation
\begin{equation}\label{resum0}
\pa(q)=C(q)+A(q)\frac{(v;q)_{\infty}}{(qu;q)_\infty}
\sum_{n=1}^\infty (-1)^n \frac{q^{2n}}{(1-2q)^n}
\frac{(uq^{n+1};q)_\infty}{(vq^n;q)_\infty},
\end{equation}
with~$A(q),C(q),u\equiv u(q),v\equiv v(q)$ rational functions given by~\eqref{acuv} and~\eqref{acuv3}.
\end{lemma}

\begin{proof}
The (easy) proof reduces to determining sufficient analyticity regions
for the various components of the basic formula~\eqref{basic},
some of the expansions being also of later use.
First, the functions $A(q)$ and $C(q)$ 
are meromorphic for~$|q|<1$, with only a pole at~$q=1/2$. They can be expanded about the point $q=1/2$ to give
\begin{align}
A & = \frac{1}{4(1-2q)^2} + \frac{1}{4(1-2q)} - \frac{1}{4}-\frac{1-2q}{4}\label{eqn:A_expansion_1/2} \\
C & = \frac{1}{4(1-2q)^2} + \frac{5}{4(1-2q)} +\frac{3}{4} - \frac{17(1-2q)}{4} + O((1-2q)^2) \label{eqn:C_expansion_1/2}
\end{align}

The function $Q(1;q)$ is analytic for $|q|<1$ except at the points for
which  $(uq;q)_{\infty}=0$, that  is,  the points  $\sigma$ for  which
$1-\sigma-\sigma^n  =0$  for $n\geq  2$.  The  smallest  of these  (in
modulus)  is  $\varphi =  (\sqrt{5}-1)/2  =  0.618034...$,  a root  of
$1-q-q^2$. So $Q(1;q)$ is  certainly analytic at $q=1/2$; the constant
term in its expansion about $q=1/2$ is
\[Q(1;1/2) = \frac{(3/2;1/2)_{\infty}}{(1/2;1/2)_{\infty}} = -0.18109782...\]

In similar fashion, $1/Q(z;q)$ is \emph{bivariate analytic} at points $(z,q)$ for which
$|q|<1$,  except when  $(vz;q)_{\infty} =  0$. This  occurs  at points
$(z_j,q)$ where  $z_j := \frac{1}{vq^j}$,
for $j\ge0$. In  particular, for $|q| < \theta$, where\footnote{
The function $v(q)=1+q^2/(1-q)$, having  nonnegative
Taylor coefficients, satisfies $|v(q)|\le v(|q|)$, for 
$|q|<1$; thus, $|1/v(q)|\ge 1/v(|q|)$. Also, $1/v(x)$ decreases from~1 to~0
for $x\in[0,1]$. Hence, with $\theta$  the real root of $1/v(\theta)=\theta$,
it follows that $|z_0|>\theta$ as soon as $|q|<\theta$.}
\begin{equation}\label{deftheta}
\theta \doteq 
0.56984 ~:=~\hbox{the  unique  real  root of  $1-2x+x^2-x^3$},
\end{equation}
we  have
$|z_0|>\theta$, hence $|z_j|>\theta$, for all~$j\ge0$.
 So,  $1/Q(z;q)$ is analytic  in the region  $\{(z,q) :
|z|,|q|  <  \theta\}$.   Thus,  for all~$n\ge 1$,
the functions  $1/Q(q^n;q)$  are  all
analytic  and uniformly bounded  by a  fixed constant,  for $|q|<r_0$,
where $r_0$ is any positive number such that $r_0<\theta$.

From these considerations, it follows that the central infinite sum that figures in~\eqref{basic}
is, when $|q|<r_1$, dominated in modulus by a positive multiple of the series
\begin{equation}\label{bound1}
\sum_n  \frac{r_1^{2n}}{(1-2r_1)^n},
\end{equation}
provided that $r_1<\theta$ \emph{and} $ r_1^2/(1-2r_1)<1$. 
Any positive $r_1$ satisfying $r_1<\sqrt{2}-1$ is then admissible.
In that case, for $|q|<r_1$, the central sum is a normally convergent sum of analytic functions; hence, 
it is analytic.
\end{proof}

The radius of analyticity of $\pa(q)$ is in fact~$1/2$. In order to obtain 
larger regions of analyticity,
one needs to improve on the reasoning underlying 
the derivation of~\eqref{bound1}.
This will result from a transformation of the central infinite sum in~\eqref{basic},
namely,
\begin{equation}\label{central}
S(q):=\sum_{n\ge1} (-1)^n \frac{q^{2n}}{(1-2q)^n}\cdot \frac{1}{Q(q^n;q)}.
\end{equation}
Only the bound $1/Q(q^n;q)=O(1)$ was used in the proof of Lemma~\ref{lem-cont0},
but we have, for instance, $1/Q(q^n;q)=1+O(q^n)$, as $n\to\infty$, and a complete expansion exists.
Indeed, since $1/Q$ is bivariate analytic in $|z|,|q|<\theta$, its $z$-expansion
at the origin is of the form 
\begin{equation}\label{qq0}
\frac{1}{Q(z;q)}=1+\sum_{\nu\ge1} d_\nu(q)z^\nu.
\end{equation}
In particular, at $z=q^n$, we have
\begin{equation}\label{qq1}
\frac{1}{Q(q^n;q)}=1+\sum_{\nu\ge1} d_\nu(q) q^{\nu n}.
\end{equation}
Now, consider the effect of an individual term $d_\nu(q)$ (instead of $1/Q(q^n;q)$)
on the sum~\eqref{central}. The identity
\begin{equation}\label{id0}
\sum_{n\ge 1} (-1)^n \frac{q^{2n}}{(1-2q)^n} q^{\nu n} =-\frac{q^{\nu +2}}{1-2q+q^{\nu +2}}
\end{equation}
provides an analytic form for the sum on the left, as long as~$q$ is not a pole of the right-hand side.
Proceeding formally, we then get,  with~\eqref{qq1} and~\eqref{id0}, upon exchanging summations
in the definition~\eqref{central} of~$S(q)$, a form of $\pa(q)$ that involves \emph{infinitely many
meromorphic elements} of the form $1/(1-2q+q^{\nu+2})$.

We shall detail validity conditions for the resulting expansion;
see~\eqref{main0} below. 
What matters, as seen from~\eqref{id0},
is the location of poles of the rational functions
$(1-2q+q^{\nu+2})^{-1}$, for $\nu\ge1$.
Define the quantities
\begin{equation}\label{defz}
\z_k := \hbox{the root 
in~$[0,1]$ of~$1-2x+x^{k+2}=0$}.
\end{equation}
We have
\[
\z_0=1; \quad 
\z_1=\frac{\sqrt{5}-1}{2}\doteq 0.618, \quad \z_2\doteq 0.543, \quad \z_3\doteq 0.518,~\ldots
\]
and~$\z_k\to\frac12$ as $k$ increases.
The location of the complex roots of $1-2x+x^{k+2}=0$ 
is discussed at length in~\cite[Ex.~V.4, p~308]{FlSe09}, as it is related to the analysis of
longest runs in binary strings: a consequence of the principle of the argument 
(or Rouch\'e's Theorem) is that,
apart from the positive real root~$\z_k$, all other complex roots lie outside the disc
$|z|<\frac34$.
The statement below builds upon this discussion and provides an extended 
analyticity region for $\pa(q)$ as well as a justification
of the validity of the expansion
resulting from~\eqref{qq1} and~\eqref{id0}, which is crucial to subsequent
developments.

\begin{lemma}\label{lem-cont1}
The generating function $\pa(q)$ is analytic at all points of the slit disc
\[
\cal D_0:=\left\{\,  q~:~
\hbox{$|q|< \frac{55}{100};~ q\not \in[\frac12,\frac{55}{100}]$}\, \right\}.
\]
For $q\in\cal D_0$, the function~$\pa(q)$ 
admits the analytic representation
\begin{equation}\label{main0}
\pa(q)=C(q)-
A(q)\frac{(v;q)_\infty}{(qu;q)_\infty}
\left[\frac{q^2}{(1-q)^2}+\sum_{\nu\ge1} d_\nu(q) \frac{q^{\nu+2}}{1-2q+q^{\nu+2}}\right],
\end{equation}
where
\[
d_\nu(q)=[z^\nu]\frac{1}{Q(z;q)}\equiv [z^\nu]\frac{(quz;q)_\infty}{(vz;q)_\infty}.
\]

In the disc $|z|<\frac{55}{100}$ punctured at~$\frac12$,
the function  $\pa(q)$ is meromorphic with simple poles at the points
$\z_2,\z_3,\ldots$, with $\z_k$ as defined in~\eqref{defz}. 
Consequently, the function $\pa(q)$ is non-holonomic,
and, in particular,  transcendental.
\end{lemma}

\begin{proof}
The starting point, noted in the proof of lemma~\ref{lem-cont0},
is that fact that $1/Q(z;q)$ is bivariate analytic 
at all points $(z,q)$ such that $|z|,|q|<\theta$, where $\theta\doteq 0.56984$ is
specified in~\eqref{deftheta}. Cauchy's coefficient formula,
\[
d_\nu(q)=\frac{1}{2i\pi} \int_{|z|=\theta_1} \frac{1}{Q(z;q)}\, \frac{dz}{z^{\nu+1}},
\]
is applicable for any $\theta_1$ such that $0<\theta_1<\theta$. 
Let us set $\theta_1=\frac{56}{100}$.
Then, since $1/Q(z;q)$ is analytic, hence continuous, 
hence bounded, for $|z|\le \theta_1$ and $|q|\le\theta_1$,
trivial bounds applied to the Cauchy integral yield
\begin{equation}\label{boundd}
|d_\nu(q)|<C\cdot \theta_1^{-\nu},
\end{equation}
for some absolute constant $C>0$. 

Consider the double sum resulting from the substitution of~\eqref{qq1} 
into~\eqref{central},
\[
S(q)=\sum_{n\ge1} (-1)^n \frac{q^{2n}}{(1-2q)^n}\cdot\bigg(1+\sum_{\nu\ge1}
d_\nu(q)q^{\nu n}\bigg).
\]
If we  constrain $q$ to be  small, say $|q|<\frac{1}{10}$, we see from~\eqref{boundd} that
the double sum is absolutely convergent. Hence, the form~\eqref{main0} is
justified for  such small values of~$q$.
We can then proceed by analytic continuation from the right-hand side
of~\eqref{main0}. The bound~\eqref{boundd} grants us the fact that 
the sum that appears there is indeed analytic in~$\cal D_0$. The 
statements, relative to the analyticity domain and the alternative
expansion~\eqref{main0} follow.
Finally, since the value~$\frac12$ corresponds to an accumulation of poles,
the function~$\pa(q)$ is non-holonomic (see, e.g.,~\cite{FlGeSa05} for context).
\end{proof}

As an immediate consequence of the dominant singularity being at~$\frac12$,
the coefficients $\pa_n$ must obey a weak asymptotic law of the form
\[
\pa_n = 2^n \theta(n), \qquad\hbox{where}\quad
\limsup_{n\to\infty}\theta(n)^{1/n}=1,\]
that is, $\theta(n)$ is a (currently unknown) subexponential factor. 

\medskip

More precise information requires a better
characterization of the behaviour of~$S(q)$, 
as $q$ approaches the dominant singularity~$\frac12$.
This itself requires a better understanding  of the coefficients $d_\nu(q)$.
To this end, we state a general and easy lemma about the
coefficients of quotients of $q$-factorials.

\begin{lemma}\label{qhyp-lemma}
 Let~$a$ be a fixed complex number satisfying $|a|<1$
and let~$q$ satisfy $|q-\frac12|<\frac1{10}$.
One has, for~$\nu\ge1$
\begin{equation}\label{qa}
[z^\nu]\frac{(az;q)_\infty}{(z;q)_\infty}=\frac{1}{(q;q)_\infty}\sum_{j=0}^\infty  
\frac{(aq^{-j};q)_\infty}{(q^{-j};q)_j}\cdot q^{j\nu}.
\end{equation}
\end{lemma}

\begin{proof}
\def\z{\overline{z}}
The function $h(z):={(az;q)_\infty}/{(z;q)_\infty}$ has simple poles at the points $\z_j:=q^{-j}$, for $j\ge0$.
We have
\[
h(z) \mathop{\sim}_{z\to\z_j} \frac{e_j(a;q)}{1-zq^j},\qquad e_j(a;q):= \frac{(aq^{-j};q)_\infty}{(q^{-j};q)_j (q;q)_\infty}.
\]
The usual expansion of coefficients of meromorphic functions~\cite[Th.~IV.10, p.~258]{FlSe09}
immediately implies a terminating form for any $J\in\Z_{\ge0}$:
\begin{equation}\label{qj}
[z^\nu]h(z)=\sum_{j=0}^J e_j(a;q) q^{j\nu}+O(R_J^n),
\end{equation}
where we may adopt $R_J=\frac32 q^{-J}$.

The last estimate~\eqref{qj} corresponds to an evaluation by residues of the Cauchy integral
representation of coefficients, 
\[
[z^\nu]h(z)=\frac{1}{2i\pi}  \, \int_{|z|=R_J} h(z)\, \frac{dz}{z^{\nu+1}}.
\]
Now, let $J$ tend to infinity. The quantity $R_J$ lies approximately midway between two
consecutive poles, $q^{-J}$ and $q^{-J-1}$, and it can be verified elementarily that,
throughout $|z|=R_J$, the function $h(z)$ remains bounded in modulus by an absolute constant
(this requires the condition~$|a|<1$). It then follows that 
we can let $J$ tend to infinity in~\eqref{qj}. For $\nu\ge1$, the coefficient integral 
taken along $|z|=R_J$ tends to~0, so that, in the limit, the exact representation~\eqref{qa} results.
\end{proof}

The formula~\eqref{qa} is equivalent to the partial fraction 
 expansion (Mittag-Leffler expansion; see~\cite[\S7.10]{Henrici74}) of the
function $h(z)$, which is meromorphic in the whole complex plane:
\begin{equation}\label{mittag}
\frac{(az;q)_\infty}{(z;q)_\infty}=1+ \frac{1}{(q;q)_\infty}\sum_{j=0}^\infty  
\frac{(aq^{-j};q)_\infty}{(q^{-j};q)_j} \frac{zq^j}{1-zq^j}.
\end{equation}
(The condition $|a|<1$ ensures the convergence of this expansion.)
As observed by Christian Krattenthaler (private communication, June 2010),
this last identity is itself alternatively deducible from the $q$-Gau{\ss} identity\footnote{ For notations, see 
Gasper and Rahman's reference text~\cite{GaRa90}: page~3 (definition of $\phi$) and  Eq.~(1.5.1), page~10
($q$-Gau{\ss} summation).}
\[
{}_2\phi_1\left[\begin{array}{c}{A,B}\\{C}\end{array};q,\frac{C}{AB}\right]=\frac{(C/A;q)_\infty (C/B;q)_\infty}
{(C;q)_\infty (C/(AB);q)_\infty},
\]
upon noticing that
\[
h(z)=\frac{(a;q)_\infty}{(1-z)(q;q)_\infty} \, {}_2\phi_1\left[\begin{array}{c} {q/a,z}\\{qz}\end{array};q,a\right].
\]

\smallskip
A direct consequence of Lemma~\ref{qhyp-lemma} is an expression
for the coefficients $d_\nu(q)=[z^\nu]Q(z;q)^{-1}$,
with $Q(z;q)$ defined by~\eqref{acuv2}:
\begin{equation}\label{dnu}
d_\nu(q)=
\frac{1}{(q;q)_\infty}\sum_{j=0}^\infty  
\frac{(quv^{-1}q^{-j};q)_\infty}{(q^{-j};q)_j}\cdot \left( vq^{j}\right)^{\nu},
\qquad \nu\ge1.
\end{equation}
To see this, set 
\[
a=quv^{-1}=\frac{q^2}{1-q+q^2},
\]
and replace $z$ by $zv$ in the definition of~$h(z)$.
Note that at~$q=1/2$, we have $u=1$, $v=3/2$, $a=1/3$, so that,
for $q\approx1/2$, we expect $d_\nu(q)$ to grow roughly 
like $(3/2)^\nu$.

Summarizing the results obtained so far, we state:

\begin{prop}\label{all1-prop}
The generating function of 3-sided prudent polygons satisfies the identity
\begin{equation}\label{all1}
\pa(q)=D(q)-q^2A(q)\frac{(a;q)_\infty (v;q)_\infty}{(q;q)_\infty (av;q)_\infty}
\sum_{\nu=1}^\infty\sum_{j=0}^\infty \left[\frac{(aq^{-j};q)_j}{(q^{-j};q)_j}
\cdot \frac{v^\nu\, q^{(j+1)\nu}}{1-2q+q^{\nu+2}}\right],
\end{equation}
where
\begin{equation}\label{defa}
a= \frac{q^2}{1-q+q^2}, \quad
v=\frac{1-q+q^2}{1-q},\quad D(q)=C(q)-\frac{q^2}{(1-q)^2}A(q)
\frac{(v;q)_\infty}{(av;q)_\infty},
\end{equation}
and $A(q), C(q)$ are rational functions 
defined in Equation~\eqref{acuv}. 
\end{prop}
\begin{proof} The identity is a direct consequence of
the formula~\eqref{dnu} for $d_\nu(q)$ and of the 
expression for~$\pa(q)$ in~\eqref{main0},
using the equivalence $av=qu$ and the simple reorganization
\[
(aq^{-j};q)_\infty =(aq^{-j};q)_j\cdot (a;q)_\infty.
\]
Previous developments imply that the identity~\eqref{all1}
 is, in particular, valid in the
real interval~$(0,\frac12)$.
The trivial equality
\begin{equation}\label{triv}
\frac{(aq^{-j};q)_j}{(q^{-j};q)_j}=
\frac{(a-q)(a-q^2)\cdots (a-q^j)}
{(1-q)(1-q^2)\cdots (1-q^j)}
\end{equation}
then shows that the expression on the right-hand side indeed represents a \emph{bona fide}
formal power series in~$q$, since the $q$-valuation of the general term of the double
sum in~\eqref{all1} increases with both~$j$ and~$\nu$. 
\end{proof}

The formula~\eqref{all1} of Proposition~\ref{all1-prop} will serve as  
the starting point of 
the asymptotic analysis of $\pa(q)$ as $q\to1/2$ in the next subsection.
Given the discussion of the analyticity of the various components in the proof of 
Lemma~\ref{lem-cont0}, the task essentially reduces to estimating the double sum in
a suitable complex neighbourhood of $q=1/2$.

\subsection{Mellin analysis}


Let $T(q)$ be the double sum that appears in the expression~\eqref{all1}
of $\pa(q)$. We shall take it here in the form
\begin{equation}\label{tq}
T=
\sum_{j=0}^\infty \frac{(aq^{-j};q)_j}{(q^{-j};q)_j} \left[
\sum_{\nu=1}^\infty 
\frac{v^\nu\, q^{(j+1)\nu}}{1-2q+q^{\nu+2}}\right]
= \sum_{j=0}^\infty \frac{(aq^{-j};q)_j}{(q^{-j};q)_j} H_j(q),
\end{equation}
with
\begin{equation}\label{hq}
H_j(q):=\sum_{\nu=1}^\infty 
\frac{v^\nu\, q^{(j+1)\nu}}{1-2q+q^{\nu+2}}.
\end{equation}
We will now study the functions~$H_j$
and propose to show that those of greater index contribute less
significant terms in the asymptotic expansion of~$\pa(q)$ near $q=1/2$.
In this way, a complete asymptotic expansion of
the function~$\pa(q)$, hence of its coefficients~$\pa_n$, can be obtained.

The main technique used here is that of Mellin transforms: we refer the reader to~\cite{FlGoDu95} for 
details of the method. The principles are recalled in~\S\ref{mel-subsubsec} below.
We then proceed to analyse the double sum~$T$ of~\eqref{tq} when $q$ is \emph{real}
and $q$ tends to $1/2$. The corresponding expansion is fairly explicit and it
is obtained at a comparatively low computational cost in~\S\ref{qreal-subsubsec}.
We finally show in~\S\ref{qall-subsubsec} that the expansion extends to a 
\emph{sector of the complex plane}
around $q=1/2$.

\def\M{\mathcal{M}}

\subsubsection{\bf\em Principles of the Mellin analysis.}\label{mel-subsubsec}
Let $f(x)$ be a complex function of the real argument~$x$. Its \emph{Mellin transform},
denoted by $f^\star(s)$ or $\M[f]$, is defined as the integral
\begin{equation}\label{meldef}
\M[f](s)\equiv f^\star(s):=\int_0^\infty f(x) x^{s-1}\, dx,
\end{equation}
where~$s$ may be complex. It is assumed that $f(x)$ is locally integrable. 
It is then well known that if $f$ satisfies the two asymptotic conditions
\[
f(x)\mathop{=}_{x\to 0} O(x^\alpha), \qquad
f(x)\mathop{=}_{x\to +\infty} O(x^\beta),
\]
with $\alpha>\beta$, then $f^\star$ is an analytic function of~$s$ in the
strip of the complex plane,
\[
-\alpha<\Re(s)<-\beta,
\]
also known as a \emph{fundamental strip}. Then, with~$c$ any real number of the interval
$(-\alpha,-\beta)$, the following inversion formula holds (see~\cite[\S VI.9]{Widder41}
for detailed statements):
\begin{equation}\label{melinv}
f(x)=\frac{1}{2i\pi}\int_{c-i\infty}^{c+i\infty} f^\star(s) x^{-s}\, ds.
\end{equation}
There are then two essential properties of Mellin transfoms.
\begin{itemize}
\item[$(\bf M_1)$] \emph{Harmonic sum property.} If the pairs $(\lambda,\mu)$ 
range over a denumerable subset of $\R\times\R_{>0}$ then one has the equality
\begin{equation}\label{harm}
\M\left[\sum_{(\lambda,\mu)} \lambda f(\mu x)\right] =
f^\star(s) \cdot \left(\sum_{(\lambda,\mu)} \lambda \mu^{-s}\right).
\end{equation}
That is to say, the harmonic sum $\sum \lambda f(\mu x)$ has a Mellin transform that 
decomposes as a product involving the transform of the base function $(f^\star)$ and 
the generalized Dirichlet series $(\sum \lambda \mu^{-s})$ associated with the ``amplitudes''~$\lambda$ and
the ``frequencies''~$\mu$. Detailed validity conditions, spelled out in~\cite{FlGoDu95},
are that the exchange of summation ($\sum$, in the definition of
the harmonic sum) and integral ($\int$, in the definition
of the Mellin transform) be permissible.
\item[$\bf (M_2)$] \emph{Mapping properties.} Poles of transforms are in correspondence with
asymptotic expansions of the original function. More precisely,
if the Mellin transform $F^\star$ of a function~$F$
admits a meromorphic extension beyond the fundamental strip, with a pole 
of some order~$m$ at some point $s_0\in\C$, with
 $\Re(s_0)<-\alpha$,
 then it contributes an asymptotic term of the form~$P(\log x)x^{-s_0}$
in  the expansion of~$F(x)$  as $x\to0$, where~$P$  is a computable polynomial of
degree~$m-1$. Schematically:
\begin{equation}\label{melasy}\def\Res{\operatorname{Res}}
F^\star(s)~\mathop{:}_{s\to s_0}~ \frac{C}{(s-s_0)^m}
\qquad\Longrightarrow\qquad F(x)~\mathop{:}_{x\to0}~ P(\log x)x^{-s_0}=\Res\left(f^\star(s)x^{-s}\right)_{s=s_0}
\end{equation}
Detailed   validity conditions,  again    spelled  out
in~\cite{FlGoDu95}, are   a     suitable decay  of     the   transform
$F^\star(s)$, as $\Im(s)\to\pm\infty$,  so as to  permit an estimate
of        the    inverse       Mellin     integral~\eqref{melinv}   by
residues -- in~\eqref{melasy}, the   expression is then none other than
the residue of $f^\star(s)x^{-s}$ at $s=s_0$.
\end{itemize}

The power of the Mellin transform  for the asymptotic analysis of sums
devolves from the   application of the mapping  property~$\bf(M_2)$ to
functions  $F(x)=\sum    \lambda    f(\mu  x)$     that  are  harmonic
sums in the sense of~$\bf(M_1)$. Indeed, the factorization property~\eqref{harm} of~$(\bf M_1)$
makes it
possible to analyse separately the singularities that arise from the base function (via $f^\star$)
and from the amplitude--frequency pairs (via $\sum \lambda\mu^{-s}$); hence
an asymptotic analysis results, thanks to $\bf (M_2)$.

\subsubsection{\bf\em Analysis for real values of~$q\to1/2$.}\label{qreal-subsubsec}
Our purpose now is to analyse the quantity~$T$ of~\eqref{tq}  with $q<1/2$,
when~$q\to1/2$. This basically reduces to analysing the quantities $H_j(q)$ of~\eqref{hq}.
Our approach consists of setting $t=1-2q$ and \emph{decoupling}\footnote{
	An instance of such a decoupling technique appears for instance
in De Bruijn's reference text~\cite[p.~27]{deBruijn81}.
} the quantities $t$ and 
$q$. Accordingly, we define the function
\begin{equation}\label{ht}
h_j(t)\equiv h_j(t;q,v):=q^{-2}\sum_{\nu=1}^\infty \frac{(vq^j)^\nu}{1+tq^{-\nu-2}},
\end{equation}
so that 
\[
H_j(q)=h_j(t;q,v(q)),\]
with the definition~\eqref{hq}. We shall let $t$ range over $\R_{\ge0}$ but restrict
the parameter~$q$ to
a small interval $(1/2-\epsilon_0,1/2+\epsilon_0)$ of $\R$ and the parameter~$v$ to 
a small interval of the form $(3/2-\epsilon_1,3/2+\epsilon_11)$, since $v(1/2)=3/2$.
We shall write such a restriction as
\[
q\approx  \frac12, \qquad  v\approx  \frac32,\] with the understanding
that $\epsilon_0,\epsilon_1$ can be taken suitably  small, as the need
arises.  Thus, \emph{for the time being, we  ignore the relations that exist
between~$t$   and  the  pair~$q,v$, and   we  shall consider   them as
independent quantities}.

As a preamble to the Mellin analysis, we state an elementary lemma

\begin{lemma} \label{strip-lem}
Let $q$ be restricted to a sufficiently small interval containing $1/2$ and~$v$ to a 
sufficiently small interval 
containing $3/2$. 
Each function $h_j(t)$ defined by~\eqref{ht} satisfies the estimate
\begin{equation}\label{defgam}
\renewcommand{\arraystretch}{1.5}
h_j(t) \mathop{=}_{t\to+\infty} O\left(\frac{1}{t}\right),\qquad 
h_j(t) \mathop{=}_{t\to0} \left\{\begin{array}{ll}
O(1) & \hbox{if $j\ge1$}\\
O(t^{-\gamma}) & \hbox{if $j=0$, with $\ds \gamma=\frac{\log v}{\log(1/q)}$}.
\end{array}\right.
\end{equation}
\end{lemma}

\noindent
For $\gamma$, we can also adopt any fixed value larger  than $\log_2(4/3)\doteq
0.415$,
provided $q$ and $v$ are taken close enough to $1/2$ and $3/2$,
respectively.

\begin{proof}
\emph{Behaviour as $t\to+\infty$.} The inequality $(1+tq^{-\nu-2})^{-1}<t^{-1}q^{\nu+2}$
implies by summation the inequality
\[
h_j(t) \le  q^{-2} t^{-1}\sum_{\nu=1}^\infty  v^\nu q^{j\nu} q^{\nu} = O\left(\frac{1}{t}\right),
\qquad t\to+\infty,
\]
given the convergence of the geometric series $\sum_\nu v q^{(j+1)\nu}$, for $v\approx3/2$
and $q\approx1/2$.

\smallskip

\emph{Behaviour as $t\to0$}.  
First, for the easy case $j\ge1$, the trivial inequality
$(1+tq^{-\nu-2})^{-1}\le 1$ implies
\[
h_j(t)=O\left(\sum_{\nu} (vq^j)^\nu\right)=O(1), \qquad t\to0.
\]
Next, for $j=0$, define the function
\[
\nu_0(t):=-2+\frac{\log (1/t)}{\log (1/q)},
\]
so that $tq^{-\nu-2}<  1$, if $\nu<\nu_0(t)$, and  $tq^{-\nu-2}\ge 1$,
if     $\nu\ge         \nu_0(t)$.      Write    $\sum_\nu            =
\sum_{\nu_0}+\sum_{\nu\ge\nu_0}$.  The sum  corresponding   to~$\nu\ge
\nu_0$ is bounded from above as in the case of $t\to+\infty$,
\[
\sum_{\nu\ge \nu_0(t)}  \frac{v^\nu}{1+tq^{-\nu-2}} \le
\sum_{\nu\ge \nu_0(t)} {v^\nu} t^{-1} q^{\nu+2} = O\left(t^{-1} (vq)^{\nu_0}\right)
 = O\left(t^{-1} (vq)^{\nu_0}\right), \qquad t\to0,\]
and the last quantity is $O(t^{-\gamma})$ for $\gamma=(\log v)/\log(1/q)$. 
The sum  corresponding to~$\nu< \nu_0$ 
is dominated by its later terms and is accordingly found to be $O(t^{-\gamma})$. 
The estimate of $h_0(t)$, as $t\to 0$, results.
\end{proof}

We can now proceed with a precise asymptotic analysis of the functions
$h_j(t)$,   as  $t\to 0$.    Lemma~\ref{strip-lem} implies  that  each
$h_j(t)$ has  its Mellin  transform  $h_j^\star(s)$ that exists  in  a
non-empty fundamental strip left   of $\Re(s)=1$. In that   strip, the
Mellin transform is
\begin{equation}\label{meli}\renewcommand{\arraycolsep}{3pt}
\begin{array}{lllll}
\M[h_j(t)] &=& 
\ds q^{-2} \M\left[\frac{1}{1+t}\right]\cdot \left(\sum_{\nu=1}^\infty (vq^j)^\nu (q^{-\nu-2})^{-s}\right) &&
\hbox{\small (by the Harmonic Sum Property ($\bf M_1$))}\\
&=&
\ds q^{-2}\M\left[\frac{1}{1+t}\right]\cdot \frac{vq^{j+3s}}{1-vq^{j+s}} &&
\hbox{\small (by summation of a geometric progression)} \\
&=&
\ds q^{-2} \frac{\pi}{\sin \pi s}  \frac{vq^{j+3s}}{1-vq^{j+s}} &&
\hbox{\small(by the classical form of $\M[(1+t)^{-1}]$).}
\end{array} 
\end{equation}
The  Mellin transform of  $(1+t)$,   which equals $\pi/\sin(\pi   s)$,
admits  $0<\Re(s)<1$  as the fundamental   strip,   so this  condition  is
necessary   for the  validity    of~\eqref{meli}.  In  addition,   the
summability of the Dirichlet  series, here plainly a geometric series,
requires the condition $|vq^{j+s}|<1$; that is,
\[
\Re(s)>-j+\frac{\log v}{\log 1/q}.
\]
In summary, the validity of~\eqref{meli} is ensured for~$s$ satisfying
\[
\lambda<\Re(s)<1, \qquad\hbox{with}\quad
\lambda:=\max\left(0,-j+\frac{\log v}{\log 1/q}\right).
\]

\begin{lemma} \label{h-lem}
For $q\approx1/2$ and $v\approx3/2$ restricted as in Lemma~\ref{strip-lem},
the function $h_j(t)$ admits an \emph{exact representation}, valid for any~$t\in(0,q^{-3})$,
\begin{equation}\label{hj}
h_j(t)=
(-1)^j\frac{vq^{3\gamma-2j-2}}{\log 1/q}
t^{j-\gamma}\Pi(\log_{1/q} t)+q^{-2}\sum_{r\ge0} (-1)^r\frac{vq^{j-3r}}{1-vq^{j-r}}t^r.
\end{equation}
Here, 
\[
\gamma \equiv \gamma(q):=\frac{\log v}{\log 1/q}
\]
so that $\gamma\approx \log_2\frac32\doteq 0.415$, when $q\approx\frac12$; the quantity~$\Pi(u)$ is an
absolutely convergent  Fourier series,
\begin{equation}\label{pi0}
\Pi(u):=\sum_{k\in\Z} p_{k} e^{-2ik\pi u},
\end{equation}
with coefficients $p_{k}$ given explicitly by 
\begin{equation}\label{pi1}
p_k=\frac{\pi}{\sin(\pi\gamma+2ik\pi^2/(\log1/q))}\,.
\end{equation}
\end{lemma}

\noindent
Observe that the $p_k$ decrease geometrically with~$k$. For instance, at~$q=1/2$, one has
\begin{equation}\label{decay}
p_k=O\left(e^{-2k\pi^2/\log 2 }\right) \doteq O\left(4.28\cdot 10^{-13}\right)^k,
\end{equation}
as is apparent from the exponential form of the sine function. Consequently,
even the very first coefficients are small: 
 at $q=1/2$, typically,
\[
|p_1|=|p_{-1}|\doteq 2.69\cdot 10^{-12}, \quad
|p_2|=|p_{-2}|\doteq 1.15\cdot 10^{-24}, \quad
|p_3|=|p_{-3}|\doteq 4.95\cdot 10^{-37}.
\]

\begin{proof}
We first perform an \emph{asymptotic} analysis of $h_j(t)$ as $t\to0^+$. This requires 
the determination of poles to the left of the fundamental strip of $h_j^\star(s)$,
and these arise from two sources.
\begin{itemize}
\item[---]
The relevant poles of $\pi/\sin\pi s$ are at $s=0,-1,-2,\ldots$; they are simple and the residue at 
$s=-r$ is $(-1)^r$.
\item[---] The quantity $(1-vq^{j+s})^{-1}$  has a simple pole
at the real point
\begin{equation}\label{sig0}
\sigma_0:=-j+\frac{\log v}{\log 1/q},
\end{equation}
as well as \emph{complex poles} of real part~$\sigma_0$,
due to the complex periodicity of the exponential function ($e^{t+2i\pi}=e^t$).
The set of all poles of $(1-vq^{j+s})^{-1}$ is then
\[
\left\{\sigma_0+\frac{2ik\pi}{\log 1/q}, \quad k\in\Z\right\}.
\]
\end{itemize}

The  proof of an  asymptotic  representation (that is, of~\eqref{hj},
with   '$\sim$' replacing  the  equality sign   there) is  classically
obtained by  integrating  $h_j^\star(s)t^{-s}$ along a  long rectangle
with corners    at $-d-iT$  and $c+iT$,   where   $c$ lies within  the
fundamental strip (in particular, between~$0$ and~$1$) and $d$ will be
taken to   be of the form   $-m-\frac12$, with $m\in  \Z_{\ge 0}$, and
smaller than $-j+\gamma$.   In   the case considered here,  there  are
regularly spaced poles along  $\Re(s)=-j+\gamma$,  so that one  should
take values  of  $T$ that are  such  that  the  line $\Im(s)=T$ passes
half-way between poles.  This, given the fast decay of $\pi/\sin\pi s$
as $|\Im(s)|$  increases and the boundedness   of the Dirichlet series
$(1-vq^{j+s})^{-1}$ along $\Im(s)=\pm T$, allows  us to let $T$ tend to
infinity.  By the  Residue   Theorem  applied to the   inverse  Mellin
integral~\eqref{melinv}, we  collect in this   way the contribution of
\emph{all} the poles at $-j+\gamma+2ik\pi/(\log 1/q)$, with $k\in\Z$, as well
as  the $m+1$  initial terms  of the  sum  $\sum_r$ in~\eqref{hj}. The
resulting  expansion  is of   type~\eqref{hj}   with the  sum $\sum_r$
truncated to $m+1$ terms and an error term that is $O(t^{m+1/2})$.

In
general,  what the   Mellin    transform    method  gives  is       an
\emph{asymptotic} rather than  \emph{exact} representation of this type.  
Here, we have more. We can finally let $m$  tend to infinity  and verify that the inverse
Mellin   integral~\eqref{melinv}  taken    along the   vertical   line
$\Re(s)=-m-\frac12$ remains  uniformly bounded in modulus  by a quantity of
the form $ct^m q^{-3m}$, for some~$c>0$.   In the limit $m\to+\infty$, the integral vanishes (as long as
$tq^{-3}<1$), and 
the exact representation~\eqref{hj} is obtained.
\end{proof}

We can now combine the \emph{identity} provided by Lemma~\ref{h-lem} with the decomposition
of the generating function~$\pa(q)$ as allowed by Equations~\eqref{tq} and  \eqref{hq},
which flow from Proposition~\ref{all1-prop}. We recall that
$H_j(q)=h_j(t;q,v(q))$.

\begin{prop}\label{all2-prop}
The generating function $\pa(q)$ of prudent polygons  satisfies,
for~$q$ in a small enough interval\footnote{
Numerical experiments suggest that in fact the formula~\eqref{mainsing}
remains valid for all $q\in(0,1/2)$.}
of the form $(1/2-\epsilon,1/2)$ (for some $\epsilon>0$),
the \emph{identity}
\begin{equation}\label{ident2}
\pa(q)=D(q)-q^2A(q)\frac{(a;q)_\infty (v;q)_\infty}{(q;q)_\infty (av;q)_\infty} T(q)
,
\end{equation}
where  the notations are those of Proposition~\ref{all1-prop},
and the function $T(q)$ admits the \emph{exact representation}
\begin{equation}\label{mainsing}
T(q)=(1-2q)^{-\gamma}\cdot \Pi\left(
\frac{\log(1-2q)}{\log 1/q}\right) U(q)+V(q),\qquad \gamma\equiv \frac{\log v}{\log 1/q},
\end{equation}
with~$\Pi(u)$ given by Lemma~\ref{h-lem}, Equations~\eqref{pi0} and~\eqref{pi1}. Set
\[
t=1-2q.
\]
The \emph{``singular series''}~$U(q)$ is
\begin{equation}\label{defU}
U(q)=\frac{vq^{3\gamma-2}}{\log 1/q}
\frac{(-q^{-1}t;q)_\infty}{(-aq^{-2}t;q)_\infty},\qquad
\gamma=\frac{\log v}{\log 1/q};
\end{equation}
and the \emph{``regular series''}~$V(q)$ is
\begin{equation}\label{defV}
V(q)=-
\frac{(q;q)_\infty}{(a;q)_\infty}\frac{q^{-2}}{1+q^{-2}t}+
q^{-2}\frac{(q;q)_\infty(av;q)_\infty}{(a;q)_\infty (v;q)_\infty}
\sum_{r=0}^\infty \frac{(a^{-1}v^{-1}q;q)_r}{(v^{-1}q;q)_r} \left(-aq^{-2}t\right)^r.
\end{equation}
\end{prop}
\begin{proof}
We start from $T(q)$ as defined by ~\eqref{tq}.
The $q$-binomial theorem is the identity~\cite[\S1.3]{GaRa90}
\begin{equation}\label{qbt}
\frac{(\theta z;q)_\infty}{(z;q)_\infty}=\sum_{n=0}^\infty \frac{(\theta;q)_n}{(q;q)_n}
z^n.
\end{equation}
Now consider the      first   term in  the    expansion~\eqref{hj}   of
Lemma~\ref{h-lem}. Sum   the corresponding  contributions for  all
values of  $j\ge0$,      after multiplication  by     the  coefficient
$\frac{(aq^{-j};q)_j}{(q^{-j};q)_j}$,  in    accordance    with~\eqref{tq}
and~\eqref{hq}. This gives
\[
U(q)=\frac{vq^{3\gamma-2}}{\log 1/q}
\sum_{j=0}^\infty \frac{(aq^{-j};q)_j}{(q^{-j};q)_j} (-q^2t)^j.
\]
A simple transformation of type~\eqref{triv} finally yields
\[
\sum_{j=0}^\infty \frac{(aq^{-j};q)_j}{(q^{-j};q)_j} (-q^{-2}t)^j
= \sum_{j=0}^\infty \frac{(a^{-1}q;q)_j}{(q;q)_j} (-aq^{-2}t)^j
\]
which provides the expression for~$U(q)$ of the singular
 series, via the $q$-binomial theorem~\eqref{qbt}
taken with $z=-at$ and $\theta=a^{-1}q$.

Summing over~$j$ in the second term in the identity~~\eqref{hj}
of
Lemma~\ref{h-lem}, we have
\[
V(q)=q^{-2}\sum_{r=0}^\infty (-q^{-2}t)^r \sum_{j=0}^\infty \frac{(aq^{-j};q)_j}{(q^{-j};q)_j}
\frac{vq^{j-r}}{1-vq^{j-r}}.
\]
Now, the Mittag-Leffler expansion~\eqref{mittag} associated with Lemma~\ref{qhyp-lemma} 
can be put in  the form
\[
\frac{(az;q)_\infty}{(z;q)_\infty}
=1+\frac{(a;q)_\infty}{(q;q)_\infty}
\sum_{j=0}^\infty \frac{(aq^{-j};q)_j}{(q^{-j};q)_j}
\frac{zq^j}{1-zq^j}.\]
An application of this identity to~$V(q)$, with $z=vq^{-r}$, shows that
\[
V(q)=q^{-2}\frac{(q;q)_\infty}{(a;q)_\infty}
\sum_{r=0}^\infty (-q^{-2}t)^r\left(\frac{(avq^{-r};q)_\infty}
{(vq^{-r};q)_\infty}-1\right),
\] 
which is equivalent to the stated form of~$V(q)$.
Note that this last form is a $q$-hypergeometric function of type ${}_2\phi_1$;
see~\cite{GaRa90}.

So far,  we have proceeded formally and  left  aside considerations of
convergence.  It  can be easily verified that  all the sums, single or
double,  involved  in  the   calculations  above are  absolutely  (and
uniformly) convergent, provided $t$  is taken small enough (i.e.,  $q$
is   sufficienty  close to~$1/2$),   given   that  all  the   involved
parameters,  such as $a,u,v$, then  stay in suitably bounded intervals
of the real line.
\end{proof}

\subsubsection{\bf\em Analysis for complex values of~$q\to1/2$.}\label{qall-subsubsec}

We  now  propose  to  show  that   the ``transcendental''   expression
of~$\pa(q)$ provided by  Proposition~\ref{all2-prop} is actually valid
in certain regions of the complex plane that extend beyond an interval
of the real line.  The regions to  be  considered are dictated  by the
requirements of the singularity analysis  method to be deployed in the
next subsection.

\begin{definition} \label{sec-def}\em
 Let $\theta_0$ be a number in the interval $(0,\pi/2)$,
called the \emph{angle},
and $r_0$ a number in~$\R_{>0}$, called the \emph{radius}.
A \emph{sector} (anchored at $1/2$) is comprised
of the set of all complex numbers $z=1/2+re^{i\theta}$ such that
\[
0<r<r_0 \qquad\hbox{and}\qquad 
\theta_0<\theta<2\pi-\theta_0.
\]
\end{definition}
We stress the   fact that the angle should   be strictly smaller  than
$\pi/2$, so that  a  sector in the   sense  of the  definition  always
includes  a part of the line  $\Re(s)=1/2$. The  smallness of a sector
will be measured by the smallness of $r_0$. That is to say:
 
\begin{prop} \label{continu-prop}
There   exists a     sector~$\cal   S_0$  (anchored   at~$1/2$),    of
angle\footnote{A careful examination of the proof 
of Proposition~\ref{continu-prop} shows that \emph{any} angle
$\theta_0>0$, however small, is suitable, but only the 
existence of \emph{some} $\theta_0<\pi/2$ is needed for singularity analysis.}   
$\theta_0<\pi/2$ and radius  $r_0>0$, such
that      the    identity    expressed by     Equations~\eqref{ident2}
and~\eqref{mainsing} holds for all $q\in\cal S_0$.
\end{prop}
\begin{proof}
The proof is a simple consequence of  analytic continuation.  We first
observe that an infinite product such as $(c;q)_\infty$ is an analytic
function   of both~$c$   and~$q$,  for   arbitrary~$c$   and  $|q|<1$.
Similarly,   the    inverse  $1/(c;q)_\infty$   is  analytic  provided
$cq^j\not=1$, for all~$c$. For instance, taking~$c=a$ where $a=a(q)=q^2/(1-q+q^2)$
and noting that $a(1/2)=1/3$, we see that $1/(a;q)$ is an analytic 
function of~$q$ in  a small complex neighbourhood of $q=1/2$.
This reasoning can  be applied to the  various Pochhammer symbols that
appear in    the   definition of   $T(q),U(q),V(q)$.    Similarly, the
hypergeometric sum  that appears in  the regular series~$V(q)$ is seen
to be analytic in the three quantities $a\approx 1/3$, $v\approx 3/2$,
and $t=1-2q\approx0$.  In  particular, \emph{the functions $U(q)$  and
$V(q)$ are analytic in a complex neighbourhood of $q=1/2$}.

Next, consider the quantity
\[
(1-2q)^{-\gamma} =\exp\left(-\gamma \log(1-2q)\right).
\]
The function~$\gamma\equiv\gamma(q)$ is  analytic  in  a neighbourhood
of~$q=1/2$,   since it  equals $(\log v)/(\log   1/q)$. The logarithm,
$\log(1-2q)$,  is analytic in any  \emph{sector} anchored at~$1/2$. By
composition, there results  that  $(1-2q)^{-\gamma}$ is  analytic in  a
small  sector anchored at~$1/2$.  It  only  remains to consider the
$\Pi$   factor   in~\eqref{ident2}.      A single    Fourier  element,
$p_ke^{-2ik\pi    u}$, with  $u=\log_{1/q}t$   and  $t=1-2q$,  is also
analytic in  a small sector (anchored at~$1/2$),  as can  be seen from
the expression
\begin{equation}\label{pkek}
p_ke^{-2ik\pi   u}=p_k      \exp\left(-2ik\pi   \frac{\log(1-2q)}{\log
1/q}\right).
\end{equation}  
Note  that, although $\Re(\log(1-2q))\to\infty$     as
$q\to1/2$, the complex exponential  
$\exp(2ik\pi\log_2(1-2q))$ remains uniformly bounded,
since $\Im(\log(1-2q))$ is bounded for~$q$ in a sector. Then, given the
fast geometric decay  of  the coefficients $p_k$ at $q=1/2$ 
(namely, $p_k=O(e^{-2k\pi^2/\log2})$; cf.~\eqref{pi1}), it follows
that  $\Pi(\log_{2}t)$  is  also  analytic  in   a sector.   
A crude adjustment of this argument (see~\eqref{sing31} and~\eqref{sing32}
below for related expansions)
suffices to verify that the geometric decay  of the terms composing~\eqref{pkek}
persists in a sector anchored at~$1/2$, so that $\Pi(\log_{1/q}t)$ is also analytic in such a sector. 

Finally, the auxiliary quantities $D(q),A(q)$ are meromorphic at~$q=1/2$,
with at most a double pole there; in particular, they are analytic in a small
enough sector anchored at~$1/2$.  
We can then choose for~$\cal S_0$ 
a  small sector that satisfies  this as well as all the  previous
analyticity constraints.
Then, by \emph{unicity of
analytic continuation}, the expression on the right-hand side of~\eqref{ident2},
with~$T(q)$ as given by~\eqref{mainsing},
\emph{must} coincide with (the analytic continuation of)
$\pa(q)$ in the sector~$\cal S_0$.
\end{proof}

\subsection{Singularity analysis and transfer}
\def\hat{\widehat}

If we drastically reduce all the non-singular quantities that occur in the main 
form~\eqref{ident2} of Proposition~\ref{all2-prop} by letting $q\to1/2$, we
are led to infer that $\pa(q)$ satisfies, in a sector
around $q=1/2$, an estimate of the form 
\begin{equation}\label{simp0}
\pa(q)= \xi_0\, (1-2q)^{-\gamma_0-2}\, \Pi(\log _2(1-2q))+O\left((1-2q)^{-3/2}\right),
\qquad \gamma_0:=\log_2(3/2),
\end{equation}
where
\begin{equation}\label{simp1}
\xi_0=-\frac1{16} U(1/2) \frac{(1/3;1/2)_\infty(3/2;1/2)_\infty}
{(1/2;1/2)_\infty (1/2;1/2)_\infty},\qquad U(1/2)=\frac{16}{9\log 2},
\end{equation}
and $U(q)$ is  the  singular series  of~\eqref{defU}.  Let us ignore for the
moment the oscillating  terms  and simplify  $\Pi(u)$ to  its constant
term $p_0$, with $p_k$ given by~\eqref{pi1}. This provides a numerical
approximation~$\hat \pa(q)$ of  $\pa(q)$. With the  general asymptotic
approximation (derived from Stirling's formula)
\begin{equation}\label{xfer0}
[q^n](1-2q)^{-\lambda} \mathop{\sim}_{n\to+\infty} 
\frac{1}{\Gamma(\lambda)} 2^n n^{\lambda-1},\qquad
\lambda\not\in\Z_{\le0},
\end{equation}
it is easily seen that $[q^n]\hat\pa(q)$ is asymptotic to the quantity
$\kappa_0 2^nn^g$ of Equation~\eqref{pnasy} in Theorem~\ref{thm-asy0},
which is  indeed the  ``principal'' asymptotic   term of $\pa_n=[q^n]\pa(q)$,
where $g=\gamma_0+1=\log_23$.

A rigorous justification and a complete analysis depend on the general
singularity  analysis     theory~\cite[Ch.~VI]{FlSe09}  applied to the
expansion  of $\pa(q)$ near~$q=1/2$.  We recall that a \emph{$\Delta$-domain
with  base~$1$} is  defined to be  the intersection  of a disc  of radius
strictly larger than~$1$ and of the complement of a sector of the form
$-\theta_0<\arg(z-1)<\theta_0$   for some $\theta_0\in(0,\pi/2)$.     A
$\Delta$-domain  with base  $\rho$  is obtained  from a $\Delta$-domain
with base~1  by means of  the  homothetic transformation $z\mapsto\rho
z$. Singularity analysis theory is then based on two types of results.
\begin{itemize} 
\item[$(\bf S_1)$] 
Coefficients   of functions in  a   basic asymptotic scale have  known
asymptotic  expansions~\cite[Th.~VI.1, p.~381]{FlSe09}. In the case of
the scale    $(1-z)^{-\lambda}$,   the  expansion, which
extends~\eqref{xfer0}, is of the form
\[
[z^n](1-z)^{-\lambda} \mathop{\sim}_{n\to+\infty}
 n^{\lambda-1}\left(1+\sum_{k\ge1} \frac{e_k}{n^k}\right),
\qquad \lambda\in\C\setminus \Z_{\ge0},
\]
where  $e_k$  is a computable   polynomial in~$\lambda$ of degree~$2k$.
Observe that this expansion is  valid for \emph{complex} values of the
exponent~$\lambda$, and if $\lambda=\sigma+i\tau$, then
\[
n^{\lambda-1} = n^{\sigma-1}\cdot n^{i\tau}=n^{\sigma-1} e^{i\tau\log n}.
\]
Thus, the real part ($\sigma$)   of the singular exponent drives   the
asymptotic  regime; the  imaginary part,  as  soon  as it  is nonzero,
induces \emph{periodic oscillations}   in the  scale of~$\log   n$.  A
noteworthy feature is that  smaller functions at the singularity $z=1$
have asymptotically smaller coefficients.

\item[$(\bf S_2)$] 
An approximation of a function near its singularity can be
transferred to an approximation of coefficients according to the rule
\[
f(z)\mathop{=}_{z\to1} O\left((1-z)^{-\lambda}\right)
\qquad\Longrightarrow\qquad
[z^n]f(z)\mathop{=}_{n\to+\infty}O\left(n^{\lambda-1}\right).
\]
The condition is that $f(z)$ be analytic in a $\Delta$--domain and 
that the $O$--approximation holds in such a $\Delta$-domain, as $z\to1$;
see~\cite[Th.~VI.3, p.~390]{FlSe09}.
Once more, smaller error terms are associated with smaller coefficients.
\end{itemize}
Equipped with  these principles, it is  possible to obtain  a complete
asymptotic expansion of   $[q^n]\pa(q)$ once a  complete  expansion of
$\pa(q)$ in the vicinity of $q=1/2$ has been obtained (set $q=z/2$, so
that  $z\approx    1$ corresponds   to~$q=1/2$).   In  this   context,
Proposition~\ref{all1-prop}  precisely grants  us the  analytic continuation of
$\pa(q)$ in a $\Delta$-domain anchored at~$1/2$, with any opening angle 
arbitrarily small; Proposition~\ref{all2-prop},
together with Proposition~\ref{continu-prop}, describe in a precise manner the
asymptotic form of $\pa(q)$ as $q\to1/2$ in a $\Delta$-domain and it is a
formal exercise to transform them into a standard asymptotic expansion,
in the form required by singularity analysis theory.

\begin{prop} \label{asy-prop}
As $q\to1/2$ in a $\Delta$-domain, the function $\pa(q)$ satisfies the
expansion
\begin{equation}\label{mainsing2}
\pa(q) \mathop{\sim}_{q\to1/2}\frac{1}{(1-2q)}+R(q)+
\sum_{j\ge1} (1-2q)^{-\gamma_0-2+j}\sum_{\ell=0}^j
\left(\log(1-2q)\right)^\ell \Pi^{(j,\ell)}(\log_2(1-2q)), \quad \gamma_0=\log_2\frac32.
\end{equation}
Here $R(q)$ is analytic at~$q=1/2$ and each $\Pi^{(j,\ell)}(u)$ is a Fourier series
\[
\Pi^{(j,\ell)}(u):=\sum_{k\in \Z} p_k^{(j,\ell)} e^{2ik\pi u},
\]
with a computable sequence of coefficients $p_k^{(j,\ell)}$.
\end{prop}
\begin{proof}
From Proposition~\ref{all2-prop}, we have
\begin{equation}\label{mainasy}
\pa(q)=\pa^{\operatorname{reg}}(q)+\pa^{\operatorname{sing}}(q),
\end{equation}
where the two terms correspond, respectively, to the ``regular'' part 
(involving the regular series~$V(q)$) and the
``singular part'' (involving the singular series~$U(q)$ as well as the factor
$(1-2q)^{-\gamma}$ and the oscillating series).

Regarding the regular part,
we have, with the notations of Proposition~\ref{all2-prop},
\begin{equation}\label{sing0}
\pa^{\operatorname{reg}}(q)=
D(q)-q^2A(q)\frac{(a;q)_\infty(v;q)_\infty}{(q;q)_\infty(q;q)_\infty}V(q).
\end{equation}
We already know that $A(q)$ and $D(q)$ are meromorphic at~$q=2$ with a
double  pole, while $V(q)$ and  the Pochhammer symbols are analytic at
$q=1/2$.   Thus,   this regular  part  has at  most   a double pole at
$q=1/2$.  A  simple computation shows that  the  coefficient  
of~$(1-2q)^{-2}$ reduces algebraically trivially -- in the sense that no $q$-identity
is involved -- to~0.
Thus, the regular part involves only  a simple pole  at $q=1/2$, as is
expressed by the first  two  terms of  the  expansion~\eqref{mainasy},
where $R(q)$ is analytic at~$q=1/2$. (Note that
the coefficient of~$(1-2q)^{-1}$ is exactly~$1$,
again for trivial reasons.) 

The singular part  is more interesting and  it can be analysed by  the
method  suggested  at  the   beginning  of this   subsection. Whenever
convenient, we freely use the abbreviation $t=1-2q$. The function
$\gamma(q)=(\log v)/(\log 1/q)$ is analytic at $q=1/2$, where
\[
\begin{array}{lll}
\gamma(q)&=&\ds \log_2\frac32+2\frac{\log 3}{(\log 2)^2}(q-\tfrac12)+\cdots\\
&\doteq&0.58496+4.5732 (q-\tfrac12)+16.317 (q-\tfrac12)^2+39.982 (q-\tfrac12)^3+86.991 (q-\tfrac12)^4
+\cdots\,.
\end{array}
\]
The function $(1-2q)^{-\gamma}$ can then be expanded as
\begin{equation}\label{expo0}
\begin{array}{lll}
\ds (1-2q)^{-\gamma(q)}&=& \ds (1-2q)^{-\gamma_0}
e^{-(\gamma(q)-\gamma_0)\log t}, \qquad\hbox{with}\quad \gamma_0=\gamma(1/2)=\log_2\frac32\\
&=& \ds (1-2q)^{-\gamma_0}\left(1+\frac{\log3}{(\log 2)^2}t \log t +
t^2 P_2(\log t)+t^3P_3(\log t)+\cdots\right),
\end{array}
\end{equation}
for a computable family of polynomials $P_2,P_3,\ldots$, where $\deg P_\ell=\ell$
and $P_\ell(0)=0$
For instance, we have, with $y:=\log t$:
\[
(1-2q)^{-(\gamma-\gamma_0)\log t}
\doteq 1+2.28ty+t^2(-4.07 y+2.61y^2)+t^3(4.99y-9.32y^2+1.99y^3)+\cdots\,.
\]
The singular series $U(q)$ of~\eqref{defU}  is analytic at $q=1/2$ and
its   coefficients   can  be determined,    both numerically  and,  in
principle,   symbolically  in terms of    Pochhammer symbols and their
logarithmic  derivatives (which   lead  to $q$-analogues of   harmonic
numbers). Numerically,  they  can be estimated  to  high precision, by
bounding the infinite  sum and products to  a finite but  large value.
(The validity of the  process can be  checked empirically by increasing
the values of this threshold, the justification being  that all
involved sums and products converge geometrically fast -- we found that
replacing $+\infty$ by 100 in numerical computations  provides 
estimates that are at least correct to 25 decimal digits.)
In this way, we obtain, for instance, 
the expansion of the function $V(q)$, which is of the form ($t=1-2q$)
\begin{equation}\label{sing2}
U(q)\doteq \frac{16}{9\log2}
+9.97\, t+21.5\, t^2+35.8\,t^3+51.9\, t^4+\cdots\,.
\end{equation}

Finally, regarding $\Pi(u)$ taken at~$u=\log_{1/q}(1-2q)$,
we note that the coefficients $p_k$ of~\eqref{pi1} can be expanded
around $q=1/2$ and pose no difficulty, while the quantities $e^{2ik\pi u}$
can be expanded by a process analogous to~\eqref{expo0}.
Indeed, we have
\begin{equation}\label{sing31}
p_k\equiv p_k(q)=\frac{\pi}{\sin(\pi\gamma_0+2ik\pi^2/\log2)}
\cdot \exp\left(
1+e_1(k)t+e_2(k)t^2+\cdots\right),
\end{equation}
where the $e_k$ only grow polynomially with~$k$. Also,
at $u=\log_{1/q}(1-2q)$, one has
\begin{equation}\label{sing32}
e^{2ik\pi u}=
(1-2q)^{2ik\pi/\log 2}
\exp\left[2ik\pi\log_2t\left(g_1t+g_2t^2+\cdots\right)\right],
\end{equation}
where the coefficients $g_j$ are those of $(\log1/q)^{-1}-(\log2)^{-1}$
expanded at $q=1/2$ and expressed in terms of~$t=1-2q$.

We can now recapitulate the results  of the discussion of the singular
part: from~\eqref{sing2},  \eqref{sing31}, and~\eqref{sing32}, we find
that the terms appearing in the singular  expansion of $\pa(q)$ are of
the form, for $j=0,1,2,\ldots$,
\[
(1-2q)^{-\gamma_0-2}t^j\left(\log t\right)^{\ell} t^{2ik\pi/\log 2},
\]
with $\ell$ such that $0\le\ell\le j$ and $k\in\Z$. The terms at fixed
$j,\ell$ add up  to form  the  Fourier series $\Pi^{(j,\ell)}$,  whose
coefficients   exhibit a  fast decrease  with~$|k|$,   similar to that
encountered   in~\eqref{decay}.  Consequently,    a finite     version
of~\eqref{mainsing2} at  any   order  holds, so  that   the  statement
results.
\end{proof}

With this last proposition, we can conclude the proof of Theorem~\ref{thm-asy0}.

\begin{proof}[Proof (Theorem~\ref{thm-asy0})]
The analytic term $R(q)$ in~\eqref{mainsing2} leaves no trace
in the asymptotic form of coefficients. Thus the global contribution of
the \emph{regular} part to coefficients $\pa_n$ reduces to~$2^n$ (with coefficient~1 
and no power of~$n$), corresponding to the term~$(1-2q)^{-1}$ in~\eqref{mainsing2}.

 The transfer to coefficients   of each term  of the
\emph{singular} part of~\eqref{mainsing2} is permissible, given the principles
of singularity analysis recalled above.  Only an amended form allowing
for logarithmic factors is needed, but this  is covered by the general
theory: for  the translation  of the  coefficients of the  basic scale
$(1-z)^{-\lambda}\log^k(1-z)$,   see~\cite[p.~387]{FlSe09}.    From   a
computational  point  of   view,  one  may  conveniently 
operate~\cite[Note   VI.7,
p.~389]{FlSe09} with
\[
[z^n](1-z)^{-\lambda}\left(\log(1-z)\right)^k =
(-1)^k \frac{\partial}{\partial\lambda}\frac{\Gamma(n+\lambda)}{\Gamma(\lambda)\Gamma(n+1)},
\]
then replace the Gamma factors of large argument
by their complete Stirling expansion.

We can now complete the proof of Theorem~\ref{thm-asy0}. 
It suffices simply to retain the terms 
corresponding to~$j=0$ in~\eqref{mainsing2}, in which case
 the error term becomes of the form $O\left((1-2q)^{-\gamma_0-1}\log(1-2q)\right)$,
which corresponds to a contribution that is $O(n^{\gamma_0}\log n)=O(n^{g-1}\log n)$
for~$\pa_n$. 

Next, regarding the Fourier element of index $k=0$, the 
function-to-coefficient correspondence yields
\[
(1-2q)^{-\gamma_0-2} \qquad \Longrightarrow\qquad  \frac{n^{\gamma_0+1}}
{\Gamma(\gamma_0+2)}
\left(1+O\bigg(\frac1n\bigg)\right).
\]
Thus, the coefficient $\kappa_0$ in~\eqref{kap0} has value 
(cf~\eqref{simp0} and~\eqref{simp1})
given by 
\[
\kappa_0={\xi_0}\cdot  \left. p_0 \right|_{q=1/2} \cdot 
\frac{1}{\Gamma(\gamma_0+2)},
\qquad \gamma_0=\frac{\log 3/2}{\log 2},
\]
with~$\xi_0$ as in~\eqref{simp1}.
This, given the form~\eqref{pi1} of~$p_k$ at $k=0$,
is equivalent to the value of~$\kappa_0$ stated in Theorem~\ref{thm-asy0}
(where $g:=\gamma_0+1=\log_23$).

For a Fourier element of index~$k\in\cal Z$, we have similarly
\[
(1-2q)^{-\gamma_0-2-i\chi_k} \qquad \Longrightarrow
\qquad \frac{ n^{\gamma_0+1+i\chi_k}}{\Gamma(\gamma_0+2+i\chi_k)}
\left(1+O(1/n)\right),\qquad\hbox{where}\quad
\chi_k:=\frac{2k\pi}{\log 2}.
\]
We finally observe that 
\[
n^{\gamma_0+1+i\chi_k}= n^{\gamma_0+1}e^{i\chi_k\log n},
\]
so that all the  terms, for $k\in\Z$,
 are  of the same asymptotic order (namely, $O(n^{\gamma_0+1})$)
and their sum constitutes a Fourier series in~$\log n$. 
The Fourier coefficient $\kappa_k$ then satisfies, from the discussion above:
\[
\kappa_k =\xi_0 \cdot \left .p_k\right|_{q=1/2} \cdot 
\frac{1} {\Gamma(\gamma_0+2+i\chi_k)}.
\]
Thus finally, with $g\equiv\gamma_0+1$:
\begin{equation}\label{fourier}
\kappa_k = 
\frac{\pi}{9\log 2\, \sin(\pi g+2ik\pi^2/\log 2)\, \Gamma(g+1+2ik\pi/\log 2)}
\prod_{j=0}^\infty \frac{(1-\frac13 2^{-j})(1-\frac32 2^{-j})}{(1-\frac12 2^{-j})^2}
.
\end{equation}

\noindent
This completes the proof of Theorem~\ref{thm-asy0}. 
\end{proof}

\begin{figure}\small
\begin{center}
\Img{5}{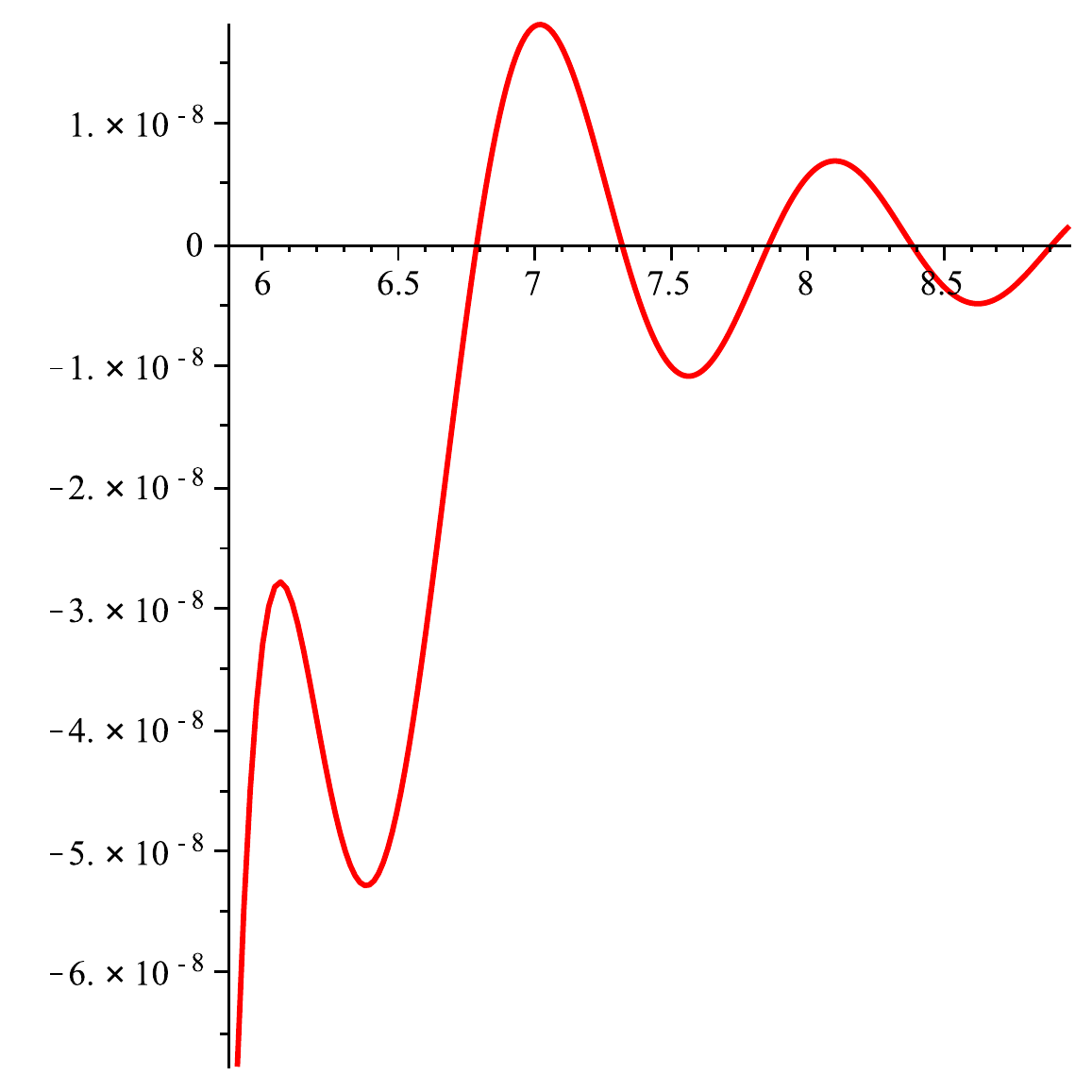}\qquad\Img{7.5}{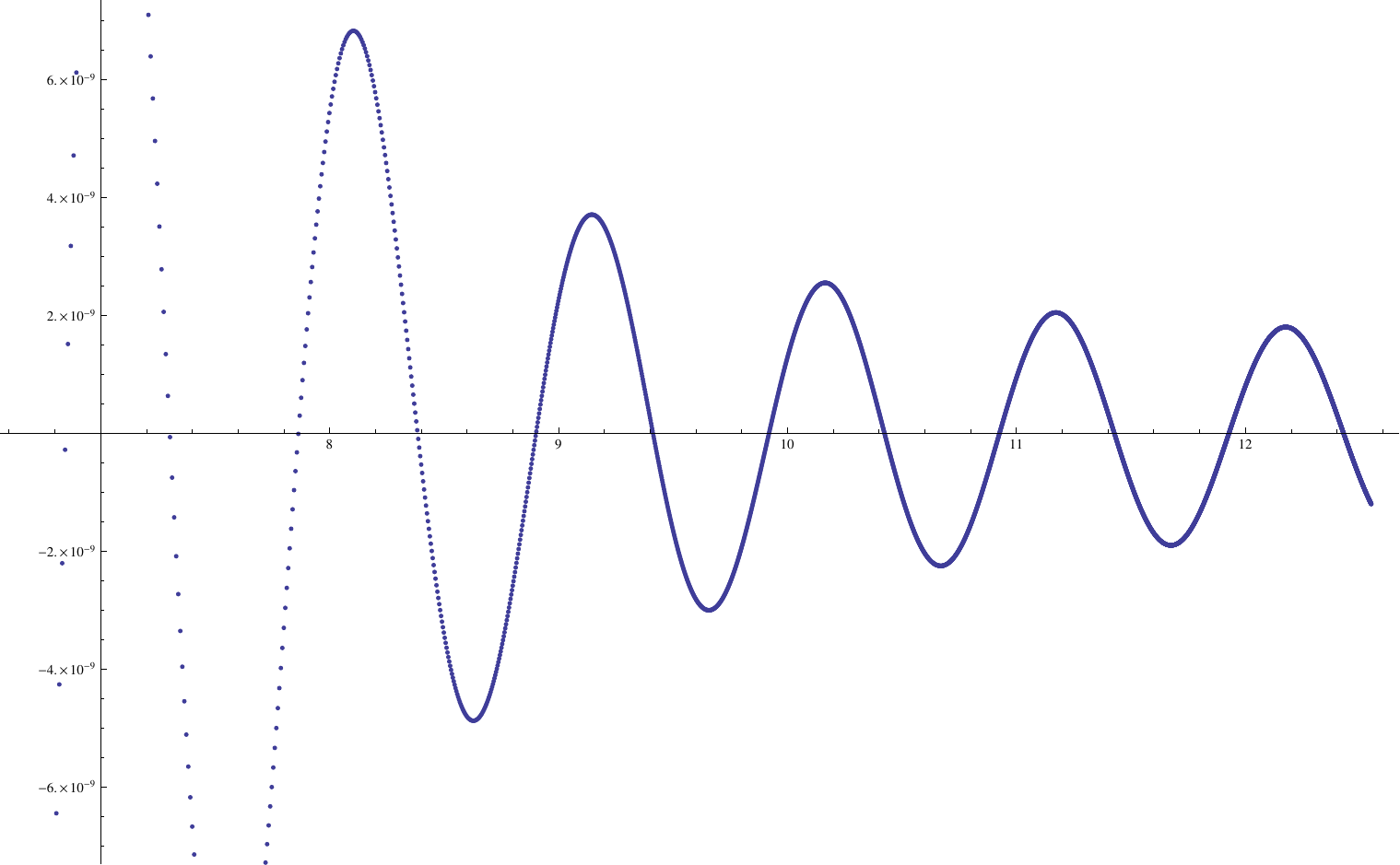}
\end{center}

\caption{\label{osci-fig}\small   
\emph{Left}: The   difference   $(\pa_n-\Omega_6)2^{-n}n^{-g}$ against $\log_2n$,
  where~$\Omega_6$ is the six-term extension of~\eqref{Omega},
  for $n=60,\ldots,500$;  the plot reveals about three periods of the
  oscillating component in the asymptotic expansion of~$\pa_n$.
\emph{Right}: The corresponding plot, for values up to $n=6000$.
The diagrams confirm the presence of
oscillations that, asymptotically,  have amplitude of the order of $ 10^{-9}$.
}
\end{figure}

The same method
shows the existence of a complete asymptotic expansion for~$\pa_n$.

\begin{theorem}\label{thm-asy1}
The number of 3-sided prudent polygons satisfies a complete asymptotic expansion,
\[
\pa_n\sim 2^n+2^n\cdot n^g\left(\Xi_{0,0}+\frac{1}{n}\left(\log n
\cdot  \Xi_{1,1}+\Xi_{1,0}\right)+\frac{1}{n^2}\left(\log^2 n
\cdot  \Xi_{2,2}+\log n
\cdot  \Xi_{2,1}+\Xi_{2,0}\right)\cdots \right),
\]
where~$\Xi_{j,\ell}$ is an absolutely convergent Fourier series in $\log n$.
\end{theorem}

The non-oscilating form obtained by retaining only the constant
terms of each Fourier series is computed by a symbolic manipulation system such
as {\sc Maple} or {\sc Mathematica} in a matter of seconds and is found to start as

\begin{small}
\renewcommand{\arraycolsep}{2pt}
\begin{equation}\label{Omega}
\begin{array}{lll}
\ds \frac{\Omega_5}{2^n}
&\doteq&  1~+~0.1083842947\,\cdot {n}^{g}
\\ &&  \ds
{}+ \left( - 0.3928066917L +
 0.5442458535 \right) \cdot n^{g-1}
\\ && 
{}+ \left(0.2627062704\, L^2
+ 0.6950193894\,L + 0.6985601031\right)\cdot n^{g-2}
\\ &&  \ds
{}+\left( 0.08310555463\, L ^{3}-
 0.02188678892\, L ^{2}-
 1.570478457\, L -1.18810811075202 \right)\cdot n^{g-3}
\\ &&  \ds
{}+\left(
 0.06722511293\, L ^{4}+ 
 0.05494834609\,L ^{3}-
 3.297513638\, L ^{2}- 4.663711650
\,L - 4.156441653\right)\cdot n^{g-4},
\end{array}
\end{equation}
\end{small}%
where $L=\log n$. In principle, all the coefficients have explicit forms in terms of the basic quantities that
appear in Theorem~\ref{thm-asy0} (augmented by derivatives of
$q$-Pochhammer symbols at small rational values).  However,
the corresponding formulae blow up exponentially, so that we only mention here the next coefficient $-0.39280\ldots$
in~\eqref{Omega}, whose exact value turns out to be
\[
-\kappa_0 \frac{\log 3}{\log^2 2}\, g.\]
Figure~\ref{osci-fig} displays the difference between $\pa_n$ and the 
six-term extension $\Omega_6$ of~\eqref{Omega}. It is piquant to note that  all the terms 
given in Equation~\eqref{Omega}  are needed in order to 
succeed in bringing the fluctuations of Figure~\ref{osci-fig} 
out of the closet. As a matter of fact, \emph{before} the analysis of Theorem~\ref{thm-asy1} was completed,
the authors had tried  empirically  to infer the likely shape
of the asymptotic expansion of~$\pa_n$, together with approximate values of coefficients,
 from a numerical analysis of series.
The success was moderate and the conclusions rather ``unstable'' (leading on occasion to heated debates between
coauthors). This state of affairs may well be present in other yet-to-be-analysed models of
combinatorics. It is pleasant that,
 thanks to complex asymptotic methods, eventually $\ldots$ everything  nicely fits in place.

\section{\bf Conclusions}\label{conclu-sec}

\emph{Classification of prudent polygons.}
Our first  conclusion is that the present  study permits us to advance
the classification of   prudent  walks and  polygons:  the  generating
functions and their coefficient asymptotics are now known in all cases
up to  3-sided (walks by length; polygons by either  perimeter or
area).  Functional equations are also known  for 4-sided prudent walks
and  polygons,   from which  it  is   possible   to  distill plausible
estimates.   We can then summarize  the present state  of knowledge by
the        following        table    (compare    with~\eqref{walk-tab}
and~\eqref{poly-tab}).

\begin{equation}\label{polya-tab}
\hbox{\small
\begin{tabular}{lcll}
\hline\hline 
 \em Polygons (area) & \em Generating function & \multicolumn{1}{c}{\em Asymptotic number} & \em References
\\
\hline
2-sided & rational & $\ds 2^{n}$ & this paper, \S\ref{2prud-subsec}
\\
3-sided & non-holonomic & $\ds C_3(n)\cdot 2^{n}n^{\log_23}$ & this paper, Th.~\ref{thm-asy0}
\\
4-sided & functional equation & 
$\ds C_4(n)\cdot 2^{n}n^{1+\log_23}~(?)$
 & this paper, \S\ref{4prud-subsec};
Beaton \emph{et al.} \cite{BeFlGu10}.
\\
\hline\hline
\end{tabular}
}
\end{equation}
(The oscillating coefficient   $C_3(n)$  is expressible in  terms  the
Fourier  series~$\kappa(u)$ of Theorem~\ref{thm-asy0}.)  The numerical
data relative  to 4-sided polygons enumerated by  area (last line) are
from~\cite{BeFlGu10}: in  all   likelihood,  the dominant  singularity
remains at $1/2$ and  the  critical exponent is~$1+\log_23$; that  is,
one more  than  the  corresponding exponent   for  3-prudent polygons.
Examination  of  subdominant terms  also suggests  that the  number of
4-polygons, once  divided by~$n$, satisfies  an expansion  of the form
obtained  in Theorem~\ref{thm-asy1}  for   the 3-prudent counterparts.
The  ``mean''  amplitude  is probably   about  0.033, and,  under  the
circumstances,  there   is little  doubt  that  minuscule oscillations
(rendered by $C_4(n)$) must  also be present.  The foregoing  analysis
of the 3-prudent case then at least has the  merit of pointing towards
the type of singularity to be expected  for 4-prudent polygons as well
as, possibly, to methods of attack for this case.

\smallskip
\emph{Methodology.}
The generating function  of 3-prudent polygons  has been found to be a
$q$-hypergeometric function, with  the argument and parameters subject
to a rational substitution. The methods  developed here should clearly
be useful in a number of similar situations.  Note that the asymptotic
enumeration of  prudent walks and polygons by  length and perimeter is
in  a way easier  since  the \emph{dominant}  singularity is polar  or
algebraic.  (Bousquet-M\'elou~\cite{Bousquet10}   however  exhibits an
interesting situation where the \emph{complete} singular structure has
a complex geometry.)

Estimates    involving periodic  oscillations  are not    unheard of in
combinatorial asymptotics~\cite{FlSe09,Odlyzko82,Odlyzko95}.  What  is
especially interesting  in the case of 3-sided prudent polygons is the pattern
of  singularities that accumulate geometrically  fast  to the dominant
singularity. This situation is prototypically encountered in
the already evoked problem of the longest run in
strings: the classical treatment is via real analytic methods followed 
by a Mellin analysis of the expressions obtained;
see~\cite{Knuth78}. 
In fact, the chain 
\begin{equation}\label{chain}
\hbox{Coefficient asymptotics $\leadsto$ 
Mellin transform + Singularity analysis}
\end{equation}
is applicable    for moment analyses. For  instance, 
the analysis of the expected longest
run of the letter~$a$ in a random binary sequence over the alphabet $\{a,b\}$
leads to the generating function~\cite[Ex.~V.4]{FlSe09}
\[
\Phi(z)=(1-z) \sum_{k\ge0} \frac{z^k}{1-2z+z^{k+1}},
\]
to which the chain~\eqref{chain} can be applied.
We could indeed recycle for our purposes some of the corresponding
analysis of~\cite{Knuth78,FlSe09}, when discussing the location of poles
of the generating function~$\pa(z)$,
in Section~\ref{sec:asymptotics}.

Another  source  of  similar phenomena   is the  analysis  of  digital
trees~\cite{Knuth98a,Szpankowski01}, when  these are   approached  via
ordinary generating functions  (rather than the customary  exponential
or Poisson generating functions).
Typically, in the simplest case of
node-depth in a random digital tree, one encounters
the generating function
\[
\Psi(z)=\frac{1}{1-z}\sum_{k\ge0} \frac{2^{-k}}{1-z(1-2^{-k})},
\]
where the geometric accumulation of poles towards~1 is transparent,
so that the chain~\eqref{chain} can once more be applied~\cite{FlRi92}.

We  should   finally  mention  that   ``critical''  exponents  similar
to~$g=\log_23$  surface  at  several  places in  mathematics,  usually
accompanied  by oscillation  phenomena,  but they  do  so for  reasons
essentially simpler than in our chain~\eqref{chain}.  For instance, in
fractal geometry, the Hausdorff dimension of the triadic Cantor set is
$1/g$, see~\cite{Falconer85}, while  that of the familiar Sierpi\'nski
gasket  is $g$,  so that~$g$  occurs as  critical exponent  in various
related     integer    sequences~\cite{FlGrKiPrTi94}.      The    same
exponent~$g=\log_23\approx 1.58$  is otherwise known  to most students
of  computer science,  since it  appears associated  to  the recurrent
sequence $f_n=n+3f_{\lfloor n/2\rfloor}$, which serves to describe the
complexity   of    Karatsuba   multiplication~\cite{Knuth98}   (where,
recursively,  the product  of \emph{two}  double-precision  numbers is
reduced to \emph{three} single-precision numbers).  In such cases, the
exponent~$g$ is eventually to be  traced to the singularity (at $s=g$,
precisely) of the Dirichlet series
\[
\omega(s)=\frac{1}{1-3\cdot 2^{-s}},
\]
which  is itself  closely  related  to the  Mellin  transforms of  our
Eq.~\eqref{meli}.     See   also   the    studies   by    Drmota   and
Szpankowski~\cite{DrSz10},      Dumas~\cite{Dumas08},      as     well
as~\cite{FlGrKiPrTi94} for elements of a general theory.

For the various reasons evoked above, we believe that the 
complex asymptotic methods developed in the present
paper are of a generality that goes somewhat beyond the 
mere case of 3-sided prudent polygons.


\def\cprime{$'$}

\end{document}